# THE HEAVY TRAFFIC LIMIT OF AN UNBALANCED GENERALIZED PROCESSOR SHARING MODEL

By Kavita Ramanan[1] and Martin I. Reiman

*Carnegie Mellon University and Alcatel–Lucent Bell Labs*

This work considers a server that processes $J$ classes using the generalized processor sharing discipline with base weight vector $\alpha = (\alpha_1, \ldots, \alpha_J)$ and redistribution weight vector $\beta = (\beta_1, \ldots, \beta_J)$. The invariant manifold $\mathcal{M}$ of the so-called fluid limit associated with this model is shown to have the form $\mathcal{M} = \{x \in \mathbb{R}_+^J : x_j = 0 \text{ for } j \in \mathcal{S}\}$, where $\mathcal{S}$ is the set of strictly subcritical classes, which is identified explicitly in terms of the vectors $\alpha$ and $\beta$ and the long-run average work arrival rates $\gamma_j$ of each class $j$. In addition, under general assumptions, it is shown that when the heavy traffic condition $\sum_{j=1}^{J} \gamma_j = \sum_{j=1}^{J} \alpha_j$ holds, the functional central limit of the scaled unfinished work process is a reflected diffusion process that lies in $\mathcal{M}$. The reflected diffusion limit is characterized by the so-called extended Skorokhod map and may fail to be a semimartingale. This generalizes earlier results obtained for the simpler, balanced case where $\gamma_j = \alpha_j$ for $j = 1, \ldots, J$, in which case $\mathcal{M} = \mathbb{R}_+^J$ and there is no state-space collapse. Standard techniques for obtaining diffusion approximations cannot be applied in the unbalanced case due to the particular structure of the GPS model. Along the way, this work also establishes a comparison principle for solutions to the extended Skorokhod map associated with this model, which may be of independent interest.

## 1. Introduction.

1.1. *Background and motivation.* Generalized Processor Sharing (GPS) is a scheduling discipline that is used to share a single processing or trans-

Received February 2006; revised February 2007.
[1]Supported in part by the NSF Grants DMS-04-06191, DMI-03-23668-0000000965, DMS-04-05343.

mission resource among traffic from several sources. Given a single server that can process one unit of work per unit of time, and that is being shared by $J$ ($1 < J < \infty$) sources or, equivalently, classes, the information needed to implement the GPS policy is contained in the weight vector $\alpha = (\alpha_1, \ldots, \alpha_J)$. When all classes have a backlog of work, class $j$ is allotted a fraction $\alpha_j$ of the total capacity of the server. When some classes achieve no backlog by using less than their allotted capacity, the remaining service capacity of the server is split among the other classes in proportion to their $\alpha_i$'s. A slight generalization of this model was considered in [17] and will be considered in this paper as well. A more precise description of that model is given here in Section 2. More background on GPS is available in [17] and the references cited there.

In [17] fluid and heavy traffic diffusion limits were obtained for the GPS model. The diffusion limits were obtained for the special "balanced case" in heavy traffic—that is, under the assumption that for every $i \in \{1, \ldots, J\}$, $\alpha_i$ is equal to the long-run average arrival rate $\gamma_i$ of work brought in by class $i$ customers. In this paper we consider diffusion limits under heavy traffic in the possibly "unbalanced case"—in which we only assume the overall condition $\sum_{i=1}^{J} \alpha_i = 1 = \sum_{i=1}^{J} \gamma_i$, without imposing any conditions on the relation between the arrival rates and weights for each class. The main appeal of the unbalanced case is the allowance of some degree of priority: classes with $\alpha_i > \gamma_i$ can be seen as receiving relatively higher priority than classes with $\alpha_i < \gamma_i$. An extreme example of this (allowed in this paper but excluded from the heavy traffic limit in [17]) is $\alpha_i = 0$, where class $i$ is served only when some other class has no backlog. The advantage of GPS with $0 < \alpha_i < \gamma_i$ is that, while class $i$ receives relatively lower priority, it cannot be completely starved, in the sense that it always receives service with rate at least $\alpha_i$. For example, as considered in [1], a "next-generation" Internet handling both real-time traffic that requires a high quality of service and best-effort traffic that has less stringent delay requirements can be modeled as a GPS system where the weight for the real-time traffic exceeds its long-run average work arrival rate, while the opposite holds for best-effort traffic.

1.2. *Relation to prior work.* In addition to the practical motivation given above, the unbalanced case is also interesting because, as elaborated below, new methods need to be developed for its analysis. There are currently two main approaches to establishing diffusion approximations for multiclass queueing networks: (i) the continuous mapping approach, which is applicable when the so-called Skorokhod map (SM), which maps the netput process to the corresponding unfinished work process, is well defined and continuous on all càdlàg paths (see, e.g., [4, 5, 20] and references therein), and (ii) the general procedure outlined in the papers [2] and [21], which is applicable when the directions of constraint (or, equivalently, directions of reflection) satisfy the so-called completely-$\mathcal{S}$ condition. Both these methods lead to diffusion approximations that are semimartingales. The present setting does not fall under either category. Indeed, neither is the GPS SM well defined for all continuous trajectories (as can be inferred from Lemma 2.4(i), Theorem 3.6 and Theorem 3.8 of [16]), nor do the GPS directions satisfy the completely-$\mathcal{S}$ condition [as can be deduced from the relation (3.1) of this paper]. In addition, the GPS diffusion limit is not in general a semimartingale.

In the balanced GPS case, this problem was circumvented in [17] by the use of the so-called extended Skorokhod map (ESM). The ESM is a generalization of the SM that allows for constraining terms that are of unbounded variation and therefore enables the pathwise construction of reflected diffusions that are not necessarily semimartingales [16]. The unfinished work process in the GPS model can be expressed as the image of the corresponding netput process under the associated GPS ESM (see Theorem 4.3 and Lemma 4.4 of [17]). Moreover, as shown in Theorem 3.6 of [16], the GPS ESM, in contrast to the GPS SM, is well defined and Lipschitz continuous on all càdlàg paths. In the balanced case, the fluid limit of the netput process is identically zero and, thus, a continuous mapping approach using the ESM, instead of the SM, can be applied to obtain diffusion approximations for the balanced GPS model (see the proof of Theorem 4.14 in [17] and also Section 5.1).

The situation in the unbalanced case is considerably more complicated because the fluid limit of the netput process is nonzero. As a result, characterization of the diffusion approximation in this case requires a better understanding of the cumulative idleness processes associated with each class. A similar generalization was considered in the context of open single-class queueing networks in [4]. The analysis in [4] hinged on certain properties of the SP associated with the queueing network studied in [4], including an explicit decomposition of the constraining term of the SP and continuity of the mapping that takes the netput process to the cumulative idleness process. As discussed in Sections 4.2 and 5.1, these properties are not satisfied in the GPS model and have consequences for both the fluid and diffusion analysis. Thus, new techniques need to be developed for the analysis of the unbalanced case.

1.3. *Main results and outline of paper.* The first step toward establishing a diffusion limit is typically the characterization of the long-time behavior of the fluid limit (see, e.g., [2] or [4]). In Section 4.2 we explicitly identify the *invariant manifold* for the fluid limit of the GPS model. The "standard definition" of an invariant manifold given, for example, in [2] (see also Definition 5.2 in [13]) is restricted to subcritical fluid limits. Our definition is slightly more general in that it also applies to supercritical fluid limits, in which case $\mathcal{M}$ provides information on how the fluid limit trajectories escape to infinity (see Remark 4.7 for further discussion of this issue). In the subcritical case, our invariant manifold can equivalently be described as the set of invariant points for the fluid limit, which coincides with the standard definition. Specifically, we show that $\mathcal{M} = \{x \in \mathbb{R}_+^J : x_i = 0 \text{ for } i \in \mathcal{S}\}$, where $\mathcal{S} \subseteq \{1, \ldots, J\}$ is the set of classes for which the long-run average work arrival rate $\gamma_i$ is strictly less than the available long-run average service rate (after reallocation of service by the GPS discipline), which we refer to as the set of strictly subcritical classes. Although it is intuitively clear that if $\gamma_i < \alpha_i$, then $i \in \mathcal{S}$, there may also exist $i \in \mathcal{S}$ for which $\gamma_i > \alpha_i$—we provide a simple, explicit characterization of the set $\mathcal{S}$. In the balanced case, this identification is trivial (there are no strictly subcritical classes), while in the unbalanced case, as mentioned above, the identification of $\mathcal{S}$ is slightly more involved in the GPS model due to the lack of an a priori representation of the constraining process in the SP decomposition, in terms of the individual cumulative idleness processes. However, we show that, when restricted to the fluid limit trajectories, it is possible to identify a certain unique decomposition of the constraining term into component processes [see Lemma 3.1 and Lemma 4.4(i)]. As shown in Theorem 4.6, this turns out to be the correct decomposition for the identification of the strictly subcritical classes and the invariant manifold $\mathcal{M}$ (see also Remark 4.5).

If the convergence of the fluid limit to $\mathcal{M}$ is sufficiently fast, one expects that $\mathcal{M}$ would also precisely characterize the state-space collapse that occurs in the heavy traffic limit. Indeed, the next main step involves proving that the strictly subcritical classes vanish in the diffusion limit, and then characterizing the behavior of the remaining, critical, classes. The lack of a simple representation for the individual cumulative idleness processes (in contrast to the networks considered in [4]) once again complicates this analysis. Nevertheless, we overcome this difficulty using two key ideas—namely, a comparison principle for the GPS ESM (Theorem 3.2), which may be of independent interest, and the introduction of the so-called reduced SP (see Theorem 5.6). This culminates in our heavy traffic limit theorem (Theorem 5.7), which is the main result of this paper. The philosophy behind our diffusion approximation is explained in greater detail in Section 5.1. We believe that our general approach is likely to be more broadly applicable to other models that do not satisfy the so-called completely-$\mathcal{S}$ condition such as, for example, networks of stations using the GPS discipline [9] or other disciplines that involve a "complete sharing" of service between classes (such as, e.g., commonly used work-conserving round-robin disciplines).

This paper is organized as follows. In the remainder of this section we describe our notational conventions. Section 2 contains a detailed description of the GPS discipline, a characterization of the unfinished work process, and definitions of the SP and ESP, culminating in a representation of the unfinished work process in terms of an SP. The fluid limit of the unbalanced GPS model is investigated and a simple characterization of the fluid limit is obtained in Section 4. Finally, Section 5 contains the statement and proof of the heavy traffic limit theorem for the unbalanced GPS model. The proofs of the results in Sections 4 and 5 rely on some basic properties of the GPS SP that are presented in Section 3.

1.4. *Notation.* We now collect together some of the notational conventions used in this paper. The set of nonnegative reals, nonnegative integers and positive integers are denoted by $\mathbb{R}_+$, $\mathbb{Z}_+$ and $\mathbb{N}$, respectively. Given $a, b \in \mathbb{R}$, $a \wedge b$ denotes the minimum of $a$ and $b$ and $a \vee b$ denotes the maximum of $a$ and $b$. Vector inequalities are to be interpreted componentwise. The standard orthonormal basis in $\mathbb{R}^J$ is represented by $\{e_i, i = 1, \ldots, J\}$, and the $J$-dimensional nonnegative orthant $\mathbb{R}^J_+$ is equal to $\{x \in \mathbb{R}^J : x \geq 0\}$. Let $\mathcal{I}$ denote the set $\{1, \ldots, J\}$. Given $E \subset \mathbb{R}^J$, $\mathcal{D}([0, \infty) : E)$ represents the space of $E$-valued right continuous functions with left limits and $\mathcal{D}_c([0, \infty); \mathbb{R}^J)$ represents the subspace of piecewise constant functions with a finite number of jumps. Unless indicated otherwise, we will assume that $\mathcal{D}([0, \infty) : E)$ and $\mathcal{D}_c([0, \infty); \mathbb{R}^J)$ are equipped with the topology of uniform convergence on compact sets (frequently abbreviated to u.o.c.). For $f \in \mathcal{D}([0, \infty) : E)$, as usual, $f(t-) = \lim_{s \uparrow t} f(s)$. For $t \in [0, \infty)$, $|f|(t)$ denotes the total variation of $f$ on $[0, t]$ with respect to the Euclidean norm $|\cdot|$ on $\mathbb{R}^J$. The composition of two functions $f$ and $g$ is as usual denoted by $f \circ g$. The identity function $\iota : [0, \infty) \to [0, \infty)$ is such that $\iota(t) = t$ for all $t \in [0, \infty)$. Given a set $A \subset \mathbb{R}^J$, $1_{\{A\}}$ represents the indicator function of the set $A$, which is defined on $\mathbb{R}^J$ and is equal to 1 on $A$ and is 0 elsewhere. In addition, $\overline{\text{co}}[A]$ denotes the closure of the convex hull of $A$, $\text{cone}[A]$ is the cone generated by $A$ and $A^\circ$ is the interior of $A$. Given a matrix $D$, we use $D'$ to denote its transpose. If $X^n, n \in \mathbb{N}$, and $X$ are processes with sample paths in $\mathcal{D}([0, \infty) : E)$, we write $X^n \Rightarrow X$ to denote weak convergence of the measures induced by $X^n$ on $\mathcal{D}([0, \infty) : E)$ to the measure induced by $X$.

**2. Model description.** In this section we provide a detailed description of the GPS discipline, introduce our assumptions on the workload arrival process, characterize the unfinished work process, and define the Skorokhod and extended Skorokhod problems.

2.1. *The GPS discipline.* We consider a single server queueing system with $J$ customer classes, where $1 < J < \infty$. Each customer arriving into the system brings in a certain amount of work that is measured in terms of the amount of time required to process it using the server's total processing capacity, which is assumed without loss of generality to be 1. The server processes the incoming work using the GPS scheduling discipline, which is described below. The work of class $i$ customers is stored in the class $i$ buffer, which is assumed to have infinite capacity. For $i \in \mathcal{I}$ and $t \in [0, \infty)$, $U_i(t)$ is defined to be the amount of work of class $i$ that is in the system at time $t$.

We now describe the GPS discipline. For $E \subseteq \mathcal{I}$, we define $\alpha_i^E$ to be the fraction of the capacity of the server that is given to class $i$ when the set of empty buffers is equal to $E$. We assume that the processor is work conserving, so that $\sum_{i \notin E} \alpha_i^E = 1$ when $E \neq \mathcal{I}$. In this paper we focus on the case when the fractions $\alpha_i^E$ are determined in the following manner by two weight vectors $\alpha \in [0,1]^J$ and $\beta \in (0,1]^J$ that satisfy $\sum_{i \in \mathcal{I}} \alpha_i = \sum_{i \in \mathcal{I}} \beta_i = 1$. Given the weight vectors, for $E = \varnothing$ (i.e., when no queue is empty), we define $\alpha_i^\varnothing \doteq \alpha_i$, and for $E \subseteq \mathcal{I}$, we let

$$\alpha_i^E \doteq \begin{cases} \alpha_i + \dfrac{\beta_i}{\sum_{j \notin E} \beta_j} \left( \sum_{j \in E} \alpha_j \right), & \text{for } i \notin E, \\ 0, & \text{otherwise.} \end{cases}$$

For all $E \subset \mathcal{I}$, $\alpha_i^E \geq \alpha_i$ for $i \in \mathcal{I} \setminus E$ and $\sum_{i \notin E} \alpha_i^E = \sum_{i \in \mathcal{I}} \alpha_i = 1$. Thus, $\alpha_i$ represents the minimum guaranteed rate assigned to class $i$. Any excess capacity is split among the remaining classes in proportion to the corresponding components of the vector $\beta$. The vectors $\alpha$ and $\beta$ will be referred to as the basic and redistribution weight vectors, respectively. The condition $\beta_i > 0$ for each $i \in \mathcal{I}$ is required to ensure that the processor is work conserving. On the other hand, we allow $\alpha_i = 0$ for some $i \in \mathcal{I}$. This represents the case when the $i$th class is of relatively low priority and only receives service when one of the high priority classes (with $\alpha_j > 0$) does not require all of its assigned capacity. More discussion on the relation between GPS and priority is contained in [17] (page 104).

2.2. *Characterization of the unfinished work process.* In this section we first introduce the primitive cumulative work arrival process associated with the GPS model and then present a characterization of the unfinished work process.

We assume that all processes are measurable functions defined on the probability space $(\Omega, \mathcal{F}, P)$. Let $H$ be the $\mathcal{D}([0,\infty) : \mathbb{R}_+^J)$-valued process such that $H_i(t)$ represents the cumulative work brought into the system by class $i$ in the interval $[0,t]$. We suppose that the cumulative work arrival process $H$ and initial conditions satisfy the following properties:

ASSUMPTION 2.1.

1. $U(0) \in \mathbb{R}_+^J$.
2. $H \in \mathcal{D}([0,\infty) : \mathbb{R}_+^J)$ is nondecreasing and piecewise constant with $H(0) = 0$.

3. $H$ has a finite number of jump points in every finite interval, almost surely.

Under this assumption, it was shown in Lemma 2.2 of [17] that the set of equations (2.1) below uniquely characterizes the set of processes $(U, I_E, E \subseteq \mathcal{I})$, where, for $E \subseteq \mathcal{I}$, $I_E(t)$ denotes the amount of time in $[0, t]$ that the set of empty buffers is equal to $E$. For $i \in \mathcal{I}$,

$$U_i(t) = U_i(0) + H_i(t) - \sum_{E \subset \mathcal{I} : i \notin E} \alpha_i^E I_E(t) \quad \text{and}$$

(2.1) $$I_E(t) = \int_0^t 1_{\{E(s) = E\}} \, ds \quad \text{for } E \subseteq \mathcal{I}, \text{ where}$$

$$E(s) \doteq \{i \in \mathcal{I} : U_i(s) = 0\}.$$

The busy time process $T_i$ defined by

(2.2) $$T_i \doteq \sum_{E \subset \mathcal{I} : i \notin E} \alpha_i^E I_E,$$

represents the cumulative amount of service given to class $i$. Note that in the definition and characterization of the unfinished work process, no assumption is made as to how the service allocated to a class is divided among the customers present in that class.

The ultimate assumptions that we require on the primitive processes are quite weak: $H$ must satisfy a functional strong law of large numbers (see Assumption 4.1) and a functional central limit theorem (see Assumption 5.1). The simplest concrete example where this happens is when the $J$ component processes of $H$ form compound renewal processes (see, e.g., [19], Lemma 2).

2.3. *Definition of the Skorokhod and extended Skorokhod problems.* The GPS model for $\alpha = \beta$ with stochastic fluid inputs was analyzed in [7]. It was shown there that the mapping taking the so-called netput process to the unfinished work can be represented in terms of a Skorokhod problem. Below, in Lemma 2.5, we recall the similar representation that was derived in [17] for the unfinished work process $U$ associated with the slight generalization of the GPS model described in Section 2.1. First, we need to recall the definitions of the Skorokhod and extended Skorokhod problems associated with the GPS model.

Roughly speaking, given the closure $G$ of a domain in $\mathbb{R}^J$, a set of allowable directions of constraint $d(x)$ associated with each point $x \in \partial G$ and a path $\psi$, the solution to the associated Skorokhod problem defines a constrained version $\varphi$ of $\psi$ that is restricted to lie in $G$ by a constraint mechanism $\varphi - \psi$ that is of bounded variation and acts in the direction of one of the vectors in $d(\varphi(s))$ using the "least effort" required to keep $\varphi$ in $G$. The solution to the extended Skorokhod problem is a generalization of the Skorokhod problem introduced in [16], which relaxes the bounded variation requirement on the constraint mechanism.

The domain of the GPS Skorokhod and extended Skorokhod problems is $G = \mathbb{R}_+^J$ and the directions of constraint (sometimes also referred to as directions of reflection) are characterized in the following manner by the redistribution weight vector $\beta \in (0,1)^J$, which satisfies $\sum_{i \in \mathcal{I}} \beta_i = 1$. Let $d_{J+1} \doteq \sum_{i \in \mathcal{I}} e_i / \sqrt{J}$ and define

$$(2.3) \qquad d_i \doteq e_i - \sum_{j \in \mathcal{I} \setminus \{i\}} \frac{\beta_j e_j}{1 - \beta_i} \qquad \text{for } i \in \mathcal{I}.$$

The set of allowable directions of constraint at any point $x$ on the boundary $\partial G$ is then given by

$$(2.4) \qquad d(x) \doteq \left\{ \sum_{i \in I(x)} a_i d_i : a_i \geq 0 \text{ for } i \in I(x) \right\},$$

where, for $x \in \mathbb{R}_+^J$, we define

$$(2.5) \qquad I(x) \doteq \begin{cases} \{i \in \mathcal{I} : x_i = 0\}, & \text{if } \sum_{j \in \mathcal{I}} x_j > 0, \\ \mathcal{I} \cup \{J+1\}, & \text{if } \sum_{j \in \mathcal{I}} x_j = 0. \end{cases}$$

The set of directions of constraint on the boundary $\partial G$ of the domain describes how service is reallocated when one or more classes have no backlog. Since this reallocation is determined solely by the redistribution weight vector $\beta$, the description of the GPS Skorokhod problem only depends on $\beta$ (and not on the basic weight vector $\alpha$). Thus, even though our model is more general than the one considered in [7, 8, 16], the results from [8, 16] can still be applied here. For more intuition on the relation between the directions of constraint and the reallocation mechanism of the GPS discipline, the reader is referred to [7, 17].

We now provide the rigorous definition of the GPS Skorokhod problem. Let $d^1(x) \doteq d(x) \cap \{x \in \mathbb{R}^J : |x| = 1\}$ and recall that for $\eta \in \mathcal{D}([0,\infty) : \mathbb{R}^J)$, $|\eta|(T)$ denotes the total variation of $\eta$ on $[0,T]$ with respect to the Euclidean norm on $\mathbb{R}^J$.

DEFINITION 2.2 (*Skorokhod problem*). Let $\psi \in \mathcal{D}([0,\infty) : \mathbb{R}^J)$ with $\psi(0) \in \mathbb{R}_+^J$ be given. Then $(\varphi, \eta)$ solves the GPS Skorokhod problem (SP) for $\psi$ if $\varphi(0) = \psi(0)$, and if for all $t \in [0, \infty)$, the following five properties hold:

1. $\varphi(t) = \psi(t) + \eta(t)$.
2. $\varphi(t) \in \mathbb{R}_+^J$.

3. $|\eta|(t) < \infty$.
4. $|\eta|(t) = \int_{[0,t]} I_{\{\varphi(s) \in \partial \mathbb{R}_+^J\}} d|\eta|(s)$.
5. There exists a measurable $\gamma \colon [0, \infty) \to \mathbb{R}^J$ such that $\gamma(t) \in d^1(\varphi(t))$ $(d|\eta|$-almost everywhere), and

$$\eta(t) = \int_{[0,t]} \gamma(s) \, d|\eta|(s).$$

Note that $\varphi$ is constrained to remain within $\mathbb{R}_+^J$, and that $\eta$ changes only when $\varphi$ is on the boundary $\partial \mathbb{R}_+^J$, in which case the change points in one of the directions in $d(\varphi)$. If $(\varphi, \varphi - \psi)$ solve the SP for $\psi$, then we denote $\varphi = \Gamma(\psi)$, and refer to $\Gamma$ as the GPS Skorokhod map (SM).

The values of $\psi \in \mathcal{D}([0, \infty) \colon \mathbb{R}^J)$ for which there exists $\varphi \in \mathcal{D}([0, \infty) \colon \mathbb{R}_+^J)$ such that $\varphi = \Gamma(\psi)$, is called the domain $\text{dom}(\Gamma)$ of the SM $\Gamma$. Theorems 1.3, 3.6 and 3.8 of [16], together, show that the domain $\text{dom}(\Gamma)$ is a strict subset of $\mathcal{D}([0, \infty) \colon \mathbb{R}^J)$ that does not include certain paths of unbounded variation. Since diffusion paths are almost surely of unbounded variation, the SM is thus inadequate for constructing reflected diffusions associated with the GPS model. This necessitates the introduction of the so-called extended Skorokhod map, first introduced in [16]. Recall that, for $A \subset \mathbb{R}^J$, $\overline{\text{co}}[A]$ represents the closure of the convex hull of the set $A$.

DEFINITION 2.3 (*Extended Skorokhod problem*). Let $\psi \in \mathcal{D}([0, \infty) \colon \mathbb{R}^J)$ with $\psi(0) \in \mathbb{R}_+^J$ be given. Then $(\varphi, \eta)$ solves the GPS extended Skorokhod problem (ESP) for $\psi$ if $\varphi(0) = \psi(0)$, and if for all $t \in [0, \infty)$:

1. $\varphi(t) = \psi(t) + \eta(t)$.
2. $\varphi(t) \in \mathbb{R}_+^J$.
3. For every $s \in [0, t]$,

(2.6) $$\eta(t) - \eta(s) \in \overline{\text{co}}\left[\bigcup_{u \in (s,t]} d(\varphi(u))\right].$$

4.

$$\eta(t) - \eta(t-) \in \overline{\text{co}}[d(\varphi(t))].$$

Theorem 3.6 of [16] shows that there exists a unique solution to the GPS ESP for all $\psi \in \mathcal{D}([0, \infty) \colon \mathbb{R}^J)$. Analogous to the GPS SM, if $(\varphi, \varphi - \psi)$ solve the GPS ESP for $\psi$, then we denote $\varphi = \overline{\Gamma}(\psi)$, and refer to $\overline{\Gamma}$ as the GPS extended Skorokhod map (ESM).

REMARK 2.4. Lipschitz continuity of the GPS SM and GPS ESM (with respect to the u.o.c. topology) on $\mathcal{D}([0, \infty) \colon \mathbb{R}^J)$ was established in Theorem 3.8 of [8] and Theorem 3.6 of [16], respectively. In particular, these results imply that solutions to the GPS SM and ESM are unique. Moreover, Theorem 1.3 of [16] shows that $\Gamma(\psi) = \overline{\Gamma}(\psi)$ for every $\psi \in \mathrm{dom}(\Gamma)$.

In order to state the Skorokhod representation for the unfinished work process, define

(2.7) $$X_i(t) \doteq U_i(0) + H_i(t) - \alpha_i t$$

and

(2.8) $$Y_i(t) \doteq \alpha_i t - T_i(t),$$

where, as in (2.2), recall that

$$T_i(t) = \sum_{E \subset \mathcal{I}: i \notin E} \alpha_i^E I_E(t).$$

Note that, by (2.1), we have $U = X + Y$.

The following result is Lemma 3.4 of [17].

LEMMA 2.5. *Let $X$ be as defined in (2.7) and let $\Gamma$ be the GPS SM associated with the weight vector $\beta \in (0,1]^J$. If $(U, I_E, E \subseteq \mathcal{I})$ satisfy (2.1), then $U = \Gamma(X)$.*

**3. Some properties of the GPS Skorokhod map.** In this section we collect properties of the GPS directions of constraint that are used for both determining the long-time behavior of the fluid limit in Section 4, as well as for identifying the diffusion limit in Section 5. A characterization of the geometry of the GPS directions of constraint is provided in Section 3.1, and a comparison principle is established in Section 3.2.

3.1. *Geometry of the directions of constraint.* The main reason that the GPS SM is not defined on all continuous functions is because the directions of constraint $d_i, i \in \mathcal{I}$, are not linearly independent (see Lemma 5.3 for one important ramification of this property and [16] for further discussion of this issue). Indeed, the directions of constraint satisfy

(3.1) $$\sum_{k=1}^{J} \beta_k (1 - \beta_k) d_k = 0,$$

which is easily verified by direct substitution. Nevertheless, as shown in Lemma 3.1, the GPS directions of constraint do exhibit a reasonable degree of regularity. This lemma is used to enable a simple description of the action of the SM on affine trajectories in Lemma 4.4, and is also used in Theorem 5.6 to establish a certain reduced representation for the GPS SM.

LEMMA 3.1 (Geometry of the GPS directions of constraint). *For every $j \in \mathcal{I}$, the vectors $\{d_1, \ldots, d_J\} \setminus \{d_j\}$ span the hyperplane $H \doteq \{x \in \mathbb{R}^J : \langle x, d_{J+1} \rangle = 0\}$. Moreover, the following two properties hold:*

(i) $H = \bigcup_{j=1}^J C_j$, where $C_j \doteq \{\sum_{k=1, k \neq j}^J \theta_k d_k : \theta_k \geq 0\}$ for $j \in \mathcal{I}$.

(ii) *Given any $w \in H$, there exists a unique vector $\theta \in \mathbb{R}_+^J$ such that $\theta_j = 0$ for at least one $j \in \mathcal{I}$ and*

$$(3.2) \qquad w = \sum_{j=1}^J \theta_j d_j.$$

*In addition, $w$ admits the representation*

$$(3.3) \qquad w_j = \frac{\theta_j}{1 - \beta_j} - \beta_j \sigma$$

*for $j \in \mathcal{I}$, where, for any set $E$ such that $\{k \in \mathcal{I} : \theta_k > 0\} \subseteq E \subset \mathcal{I}$,*

$$(3.4) \qquad 0 \leq \sigma \doteq \sum_{k:\theta_k>0} \frac{\theta_k}{1-\beta_k} = \sum_{k \in E} \frac{\theta_k}{1-\beta_k} = \frac{\sum_{j \in E} w_j}{1 - \sum_{j \in E} \beta_j}.$$

PROOF. The first statement of the lemma follows from the proof of Lemma 3.1 of [8]. In turn, this statement implies that $\bigcup_{j=1}^J C_j \subseteq H$ and that, given any $w \in H$, there exist $\theta_k \in \mathbb{R}$, $k \in \{2, \ldots, J\}$, such that $w = \sum_{k=2}^J \theta_k d_k$. If $\theta_k \geq 0$ for every $k \in \{2, \ldots, J\}$, then it automatically follows that $w \in C_1$. Otherwise, choose $l \in \{2, \ldots, J\}$ such that

$$l = \arg\min_{k \in \{2,\ldots,J\}} \frac{\theta_k}{\beta_k(1 - \beta_k)},$$

define $\theta_1 \doteq 0$ and, using relation (3.1), eliminate $d_l$ in the representation of $w$ to obtain

$$w = \sum_{k=1, k \neq l}^J \left( \theta_k - \theta_l \frac{\beta_k(1-\beta_k)}{\beta_l(1-\beta_l)} \right) d_k.$$

In this case, by the choice of $l$, $\theta_l < 0$ and all the coefficients of $d_k, k \neq l$, in the above decomposition are nonnegative. This shows that $w \in C_l$. So, in both cases, we have shown that $w \in \bigcup_{j=1}^J C_j$. Since this holds for every $w \in H$, we conclude that $H \subseteq \bigcup_{j=1}^J C_j$, which, when combined with the reverse inclusion, establishes property (i).

For property (ii), fix $w \in H$ and first observe that the existence of $\theta \in \mathbb{R}_+^J$ satisfying (3.2) follows immediately from property (i). Now suppose that there exist $i, j \in \mathcal{I}$, $\theta \in \mathbb{R}_+^J$ with $\theta_i = 0$ and $\tilde{\theta} \in \mathbb{R}_+^J$ with $\tilde{\theta}_j = 0$ such that

$$(3.5) \qquad w = \sum_{k=1, k \neq i}^J \theta_k d_k = \sum_{k=1, k \neq j}^J \tilde{\theta}_k d_k.$$

If $i = j$, then the linear independence of the vectors $d_k, k \in \mathcal{I} \setminus \{i\}$ shows that $\theta = \tilde{\theta}$. Now suppose that $i \neq j$. Then, since relation (3.1) implies that

$$d_j = \frac{-1}{\beta_j(1-\beta_j)} \sum_{k=1, k \neq j}^{J} \beta_k(1-\beta_k)d_k,$$

substituting this equality into the first representation for $w$ in (3.5) yields

$$w = \sum_{k=1, k \neq i}^{J} \theta_k d_k = \sum_{k=1, k \notin \{i,j\}}^{J} \theta_k d_k + \theta_j d_j$$

$$= \sum_{k=1, k \notin \{i,j\}}^{J} \left[\theta_k - \theta_j \frac{\beta_k(1-\beta_k)}{\beta_j(1-\beta_j)}\right] d_k - \theta_j \frac{\beta_i(1-\beta_i)}{\beta_j(1-\beta_j)} d_i.$$

Comparing the last display with the second representation for $w$ in (3.5), the linear independence of the vectors $d_k, k \in \mathcal{I} \setminus \{j\}$, implies that

$$\tilde{\theta}_k = \begin{cases} \theta_k - \theta_j \dfrac{\beta_k(1-\beta_k)}{\beta_j(1-\beta_j)}, & \text{for } k \in \mathcal{I} \setminus \{i,j\}, \\ -\theta_j \dfrac{\beta_i(1-\beta_i)}{\beta_j(1-\beta_j)}, & \text{for } k = i. \end{cases}$$

The fact that $\tilde{\theta}_i \geq 0$, $\theta_j \geq 0$ and $\beta_i(1-\beta_i)/[\beta_j(1-\beta_j)] > 0$ together imply that $\theta_j = 0$, $\tilde{\theta}_i = 0$ and, hence, that $\tilde{\theta}_k = \theta_k$ for $k \in \mathcal{I} \setminus \{i,j\}$. This establishes uniqueness of the vector $\theta$ satisfying (3.2).

To show the last assertion, using the definition of the directions of constraint and the first property of (ii), we obtain

$$w_j = \frac{\theta_j}{1-\beta_j} - \beta_j \sum_{k:\theta_k>0} \frac{\theta_k}{1-\beta_k}$$

for $j \in \mathcal{I}$, which proves (3.3) with $\sigma$ as defined in the first equality of (3.4). The second equality follows trivially since $\theta_k = 0$ for $k \in E \setminus \{k:\theta_k>0\}$. Summing the last display over all $j$ such that $\theta_j > 0$, we obtain

$$\sum_{j:\theta_j>0} w_j = \left(1 - \sum_{j:\theta_j>0} \beta_j\right)\left(\sum_{k:\theta_k>0} \frac{\theta_k}{1-\beta_k}\right),$$

which proves the third equality in (3.4) when $E = \{j:\theta_j>0\}$. Since (3.3) shows that for every $j$ with $\theta_j = 0$, $w_j = \sigma \beta_j$, this now shows that the third equality in (3.4) holds for any $E \subset \mathcal{I}$ that contains $\{j:\theta_j>0\}$. □

3.2. *A comparison principle.* The main result of this section is a comparison principle for the GPS ESM, which may be of more general interest. Analogous comparison results for the one-dimensional SM on $[0,\infty)$ and the two-sided SM on $[0,a]$ can be found in [22] and [14], respectively. The statement of the comparison principle involves the projection operator $\pi:\mathbb{R}^J \to \mathbb{R}^J_+$ associated with the GPS SP, which is characterized by the property that

(3.6)
$$\pi(x) = x \quad \text{if } x \in \mathbb{R}^J_+,$$
$$\pi(x) \in \partial \mathbb{R}^J_+ \quad \text{and} \quad \pi(x) - x \in d(\pi(x)) \quad \text{if } x \notin \mathbb{R}^J_+.$$

The existence, uniqueness and Lipschitz continuity of the GPS projection operator was established in Theorem 3.8 of [8]. The comparison principle of Theorem 3.2 is used in Theorem 4.6 to provide uniform bounds (with respect to initial conditions in a compact set) on the time at which the fluid limit of the GPS model reaches the invariant manifold. It is also used in Section 5.2 to establish "state space collapse" for the diffusion limit (see Theorem 5.6).

THEOREM 3.2 (A comparison principle). *Let $\pi$ be the GPS projection operator and $\overline{\Gamma}$ be the GPS ESM. If $\tilde{\nu} \leq \nu$, then*

(3.7) $$\pi(\tilde{\nu}) \leq \pi(\nu).$$

*Moreover, if for $\tilde{\psi}, \psi \in \mathcal{D}([0,\infty); \mathbb{R}^J)$, there exists a coordinate-wise nondecreasing function $\Delta$ with $\Delta(0) \geq 0$ such that $\psi(t) = \tilde{\psi}(t) + \Delta(t)$ for every $t \in [0,\infty)$, then*

(3.8) $$\overline{\Gamma}(\tilde{\psi})(t) \leq \overline{\Gamma}(\psi)(t) \quad \text{for every } t \in [0,\infty).$$

In the following lemma, we first show how (3.8) can be deduced from (3.7) for a general class of ESPs.

LEMMA 3.3. *Let $\pi$ be the projection operator associated with an ESP that has a uniformly continuous ESM defined on $\mathcal{D}([0,\infty); \mathbb{R}^J)$. Suppose that $\tilde{\nu} \leq \nu$ implies $\pi(\tilde{\nu}) \leq \pi(\nu)$. Then, for $\tilde{\psi}, \psi \in \mathcal{D}([0,\infty); \mathbb{R}^J)$ such that there exists a coordinate-wise nondecreasing function $\Delta$ with $\Delta(0) \geq 0$ and $\psi(t) = \tilde{\psi}(t) + \Delta(t)$ for every $t \in [0,\infty)$, (3.8) holds. Moreover, if the SM is uniformly continuous and well defined on $\mathcal{D}([0,\infty); \mathbb{R}^J)$, then (3.8) holds with the ESM replaced by the SM.*

PROOF. When $\tilde{\psi}, \psi \in \mathcal{D}_c([0,\infty); \mathbb{R}^J)$ (recall that $\mathcal{D}_c([0,\infty); \mathbb{R}^J)$ is the subspace of $\mathcal{D}([0,\infty); \mathbb{R}^J)$ that has piecewise constant trajectories with a finite number of jumps), this can be proved using induction. Indeed, if $0 =$ $t_0 < t_1 < t_2 < t_n$ are the union of the jump points of $\psi$ and $\tilde{\psi}$, it follows that $\overline{\Gamma}(\psi)(0) = \pi(\psi(0)), \overline{\Gamma}(\tilde{\psi})(0) = \pi(\tilde{\psi}(0))$ and for $k = 1, \ldots, n-1$,

$$\overline{\Gamma}(\psi)(t_{k+1}) = \pi(\overline{\Gamma}(\psi)(t_k) + \psi(t_{k+1}) - \psi(t_k)),$$

and likewise for $\overline{\Gamma}(\tilde{\psi})(t_{k+1})$ (see, e.g., (28) of [6] for this construction for the SM $\Gamma$—it is easy to see that the same construction also holds for the ESM $\overline{\Gamma}$). Now suppose $\overline{\Gamma}(\tilde{\psi})(t_k) \leq \overline{\Gamma}(\psi)(t_k)$. By the assumption on $\Delta$, we have $\tilde{\psi}(t_{k+1}) - \tilde{\psi}(t_k) \leq \psi(t_{k+1}) - \psi(t_k)$, and thus, (3.7) implies that $\overline{\Gamma}(\tilde{\psi})(t_{k+1}) \leq \overline{\Gamma}(\psi)(t_{k+1})$. Since $\Delta(0) \geq 0$, (3.7) ensures that $\overline{\Gamma}(\tilde{\psi})(t_0) \leq \overline{\Gamma}(\psi)(t_0)$ and thus, by induction, (3.8) holds whenever $\psi, \tilde{\psi} \in \mathcal{D}_c([0, \infty); \mathbb{R}^J)$.

In order to extend the result to general $\psi, \tilde{\psi}, \Delta$ in $\mathcal{D}([0, \infty); \mathbb{R}^J)$ that satisfy the assumptions of the lemma, we need only use the fact that $\psi, \tilde{\psi}, \Delta$ can each be approximated in the uniform norm by corresponding sequences $\{\psi^{(n)}\}, \{\tilde{\psi}^{(n)}\}$ and $\{\Delta^{(n)}\}$ of piecewise constant trajectories with a finite number of jumps such that $\psi^{(n)} = \tilde{\psi}^{(n)} + \Delta^{(n)}$, where $\Delta^{(n)}(0) = 0$ and $\Delta^{(n)}$ is coordinate-wise nondecreasing. Indeed, consider the sequence $l_n \doteq \{0 = t_0 < t_1 < \cdots < t_{k_n}^n, k_n \in \mathbb{N}\}$ of partitions of $[0, \infty)$, such that the $n$th partition contains all points at which $\psi$ or $\psi'$ has a jump of magnitude greater than $1/n$ and the mesh size $\max_{i=1,\ldots,k_n} |t_i^n - t_{i-1}^n| \to 0$ as $n \to \infty$. Define

$$f^{(n)}(t) = f(t_{i-1}) \qquad \text{for } t \in [t_{i-1}^n, t_i^n), \ i = 1, \ldots, k_n,$$

for $f \in \{\psi, \tilde{\psi}, \Delta\}$. Then for every $n \in \mathbb{N}$, $f^{(n)} \to f$ u.o.c. for $f \in \{\psi, \tilde{\psi}, \Delta\}$ and $\psi^{(n)} = \tilde{\psi}^{(n)} + \Delta^{(n)}$ with $\Delta^{(n)}(0) = 0$, $\Delta^{(n)} \geq 0$ and $\Delta^{(n)}$ coordinate-wise nondecreasing. Thus, the inequality (3.8) is satisfied with $\psi, \tilde{\psi}$ replaced by $\psi^{(n)}$ and $\tilde{\psi}^{(n)}$, respectively, for every $n \in \mathbb{N}$. Taking limits as $n \to \infty$, the uniform continuity of $\overline{\Gamma}$ then ensures that the inequality (3.8) also holds for $\psi, \tilde{\psi}$. If the SM is defined and uniformly continuous on $\mathcal{D}([0, \infty); \mathbb{R}^J)$, then the last statement of the lemma follows due to exactly the same argument as that used for the ESM. This completes the proof of the lemma. □

Before presenting the proof of the theorem, it will be convenient to introduce the following notation. Recall that the one-dimensional Skorokhod mapping, $\Gamma_1 : \mathcal{D}([0, \infty) : \mathbb{R}) \to \mathcal{D}([0, \infty) : \mathbb{R})$, is given explicitly by

$$(3.9) \qquad \Gamma_1(f)(t) = f(t) + \sup_{s \in [0,t]} [-f(s)] \vee 0.$$

Also, given $j \in \mathcal{I}$, define $P^{(j)}$ to be the matrix whose $i$th column is $d_i$ for $i \neq j$, and whose $j$th column is the unit vector $e_j$. Let $Q^{(j)} = I - P^{(j)}$. Then it follows immediately from the definition (2.3) of $d_i$ that $Q_{ik}^{(j)} \geq 0$ for $i, k \in \mathcal{I}$, $\sum_{i \in \mathcal{I}} Q_{ij}^{(j)} = 0$ and for every $k \in \mathcal{I} \setminus \{j\}$, $\sum_{i \in \mathcal{I}} Q_{ik}^{(j)} = 1$. We now claim that $(Q^{(j)})'$ is the transition matrix of a transient $J$-state sub-Markov chain, where the chain "dies" after entering state $j$. To see that it is transient, note that for $i \neq j, (Q^{(j)})'_{ij} = \beta_j/(1-\beta_i) > 0$ (since, by assumption, $\min_{j \in \mathcal{I}} \beta_j > 0$). Thus, $(Q^{(j)})'$ is strictly substochastic and, hence, has spectral radius $\sigma(Q^{(j)}) < 1$. We now describe the SP associated with the matrix $P^{(j)}$. Roughly speaking, it has domain $\mathbb{R}^J_+$ and direction of constraint in the relative interior of the $i$th boundary face, $\{x \in \mathbb{R}^J_+ : x_j = 0\}$, given by the $i$th column of $P^{(j)}$. More precisely, for $j \in \mathcal{I}$, define the set-valued function $d^{(j)}(\cdot)$ on $\partial \mathbb{R}^J_+$ as follows:

$$d^{(j)}(x) = \begin{cases} d(x), & \text{if } x_j > 0, \\ \left\{ a_j e_j + \sum_{i \in \mathcal{J}(x) \setminus \{j\}} a_i d_i : a_i \geq 0 \text{ for } i \in \mathcal{J}(x) \right\}, & \text{if } x_j = 0, \end{cases}$$

where, for $x \in \mathbb{R}^J_+$, $\mathcal{J}(x) \doteq \{i \in \mathcal{I} : x_i = 0\}$. Also, let $d^{(j)^1}(x) \doteq d^{(j)}(x) \cap \{v \in \mathbb{R}^J : |v| = 1\}$. The SP associated with $P^{(j)}$ is defined as in Definition 2.2, but with $d^1(\cdot)$ replaced by $d^{(j)^1}(\cdot)$. Let $\Gamma^{(j)}$ be the corresponding SM, and let $\pi^{(j)}$ be the corresponding projection operator, characterized by the relations (3.6), but with $d(\cdot)$ replaced by $d^{(j)}(\cdot)$. SPs of this kind were introduced in [11], the results of which show that the properties of $P^{(j)}$ described above guarantee that $\Gamma^{(j)}$ and $\pi^{(j)}$ are well defined on $\mathcal{D}([0,\infty);\mathbb{R}^J)$ and $\mathbb{R}^J$, respectively (see also the discussion in Section 2.3 of [8]).

We now present the proof of Theorem 3.2.

PROOF OF THEOREM 3.2. In order to prove the first property, fix $\tilde{\nu}, \nu \in \mathbb{R}^J$ such that $\tilde{\nu} \leq \nu$ and define $\tilde{\kappa} = \pi(\tilde{\nu})$ and $\kappa = \pi(\nu)$. Then the definitions of $\pi$ and the GPS directions of constraint show that

$$(3.10) \quad \sum_{i \in \mathcal{I}} \tilde{\nu}_i \leq 0 \implies -\tilde{\nu} \in d(0) \implies \pi(\tilde{\nu}) = 0.$$

Since $\pi(\nu) \geq 0$, this proves the result when $\sum_{i \in \mathcal{I}} \tilde{\nu}_i \leq 0$. We shall consider two cases with $\sum_{i \in \mathcal{I}} \tilde{\nu}_i > 0$.

*Case* 1. $\sum_{i \in \mathcal{I}} \tilde{\nu}_i > 0$ and there exists $j \in \mathcal{I}$ such that $\tilde{\kappa}_j > 0$ and $\kappa_j > 0$.

In this case, we first claim that $\pi(\nu) = \pi^{(j)}(\nu)$. Indeed, this is a direct consequence of the definition and uniqueness of the projection operators $\pi$ and $\pi^{(j)}$, and the fact that $\kappa_j > 0$ implies $d^{(j)}(\kappa) = d(\kappa)$. Since $\tilde{\kappa}_j > 0$, the same argument also shows that $\pi(\tilde{\nu}) = \pi^{(j)}(\tilde{\nu})$. From the definition of the SP and the projection operator, it is easy see that $\Gamma^{(j)}(\nu \iota) = \pi^{(j)}(\nu)\iota$. Moreover, the comparison principle proved in Theorem 4.1 of [18] guarantees that $\Gamma^{(j)}(\tilde{\nu}\iota)(t) \leq \Gamma^{(j)}(\nu\iota)(t)$ for $t \in [0,\infty)$, and so, substituting $t = 1$, we obtain the inequality $\pi^{(j)}(\tilde{\nu}) \leq \pi^{(j)}(\nu)$. Thus, $\pi(\tilde{\nu}) = \pi^{(j)}(\tilde{\nu}) \leq \pi^{(j)}(\nu) = \pi(\nu)$, and so the theorem is true in this case.

*Case* 2. $\sum_{i \in \mathcal{I}} \tilde{\nu}_i > 0$ and $\{j \in \mathcal{I} : \kappa_j = 0\} = \{j \in \mathcal{I} : \tilde{\kappa}_j > 0\}$.

We shall argue by contradiction to show that this case cannot occur when $\tilde{\nu} \leq \nu$. For $\nu \in \mathbb{R}^J$, let $\mathcal{E}(\nu) \doteq \{j \in \mathcal{I} : \kappa_j = 0\}$. Note that the first condition of Case 2 implies that $-\tilde{\nu} \notin d(0)$ and, hence, $\tilde{\kappa} \neq 0$ or, equivalently, $\mathcal{E}(\tilde{\nu}) \neq \mathcal{I}$. Next, suppose $\mathcal{E}(\nu) = \mathcal{I}$. Then we have $\kappa = 0$ and so, by the definition of $\pi$, $\sum_{i \in \mathcal{I}} \nu_i \leq 0$. When combined with the ordering $\tilde{\nu} \leq \nu$, this contradicts the first assumption of the case. Now consider the remaining possibility when $\mathcal{E}(\nu) \subset \mathcal{I}$ and $\mathcal{E}(\tilde{\nu}) \subset \mathcal{I}$. Then, by the definition of $\pi$ and Lemma 3.1, $w \doteq \kappa - \nu \in d(\kappa) = \overline{\text{cone}}[d_j, j \in \mathcal{E}(\nu)] \subseteq H$. This means, in particular, that if $\theta \in \mathbb{R}^{J+1}$ is the unique vector in (3.2), then $\theta_{J+1} = 0$ and $\{j : \theta_j > 0\} \subseteq \mathcal{E}(\nu)$. Thus, using (3.2)–(3.4) and the fact that $\nu_j = -w_j$ for $j \in \mathcal{E}(\nu)$, we observe that

$$(3.11) \qquad \sum_{j \in \mathcal{E}(\nu)} \nu_j = -\sum_{j \in \mathcal{E}(\nu)} w_j = -\left(1 - \sum_{j \in \mathcal{E}(\nu)} \beta_j\right)\sigma \leq 0.$$

Since $\mathcal{E}(\tilde{\nu}) \neq \mathcal{I}$, an analogous argument also shows that $\sum_{j \in \mathcal{E}(\tilde{\nu})} \tilde{\nu}_j \leq 0$. Moreover (3.11), along with the fact that $\tilde{\nu} \leq \nu$, implies that

$$\sum_{j \in \mathcal{E}(\nu)} \tilde{\nu}_j \leq \sum_{j \in \mathcal{E}(\nu)} \nu_j \leq 0.$$

On the other hand, since $\mathcal{E}^c(\nu) = \mathcal{E}(\tilde{\nu})$ due to the second condition of Case 2, we also have

$$\sum_{j \in \mathcal{E}^c(\nu)} \tilde{\nu}_j = \sum_{j \in \mathcal{E}(\tilde{\nu})} \tilde{\nu}_j \leq 0.$$

Together, the last two relations imply that $\sum_{j \in \mathcal{I}} \tilde{\nu}_j \leq 0$, which contradicts the first assumption of the case. This completes the proof of (3.7).

Since the GPS ESM is Lipschitz continuous, the second assertion of the theorem follows from the first due to Lemma 3.3. □

**4. Long-time behavior of the fluid limit.** In this section we consider a sequence of GPS systems with an associated sequence $\{H^n\}$ of cumulative work arrival processes defined on $(\Omega, \mathcal{F}, P)$ that satisfy Assumption 2.1. Let $U^n$ and $T^n$ be the associated unfinished work and busy time processes uniquely characterized by equations (2.1) and (2.2). We also consider the associated sequences $\{X^n\}$ and $\{Y^n\}$, where $X^n$ and $Y^n$ are defined by the relations (2.7) and (2.8), respectively, with $U_i(0)$, $H_i$ and $T_i$ replaced by $U_i^n(0)$, $H_i^n$ and $T_i^n$, respectively. Assume $\mathcal{F}$ is complete with respect to $P$ and for $n \in \mathbb{N}$, let $\{\mathcal{F}_t^n\}$ be complete filtrations such that $H^n$ is adapted to $\{\mathcal{F}_t^n\}$.

In Section 4.1 we state the characterization of the fluid limit of the sequence of unfinished work processes that was obtained in [17]. In Section 4.2 we identify the invariant manifold for the fluid limit—this constitutes the first step toward establishing the heavy-traffic diffusion approximation in Section 5.

4.1. *Characterization of the fluid limit of the unfinished work process.*
Given a sequence $\{f^n\} \subset \mathcal{D}([0,\infty):\mathbb{R}^J)$, we define the associated fluid scaled sequence $\{\overline{f}^n\} \subset \mathcal{D}([0,\infty):\mathbb{R}^J)$ by

$$\overline{f}^n(t) \doteq \frac{f^n(nt)}{n} \qquad \text{for } t \in [0,\infty). \tag{4.1}$$

From Definition 2.2 and the fact that the GPS SM $\Gamma$ is Lipschitz continuous on its domain (see Remark 2.4), it is easy to verify that $\Gamma$ is nonanticipatory in the sense that $\Gamma(X)(t)$ depends only on $\{X(s), s \leq t\}$. Since $H^n$, and therefore $X^n$, is adapted to the filtration $\{\mathcal{F}^n_t\}$ and since $U^n$ is right continuous and $U^n = \Gamma(X^n)$ by Lemma 3.4 of [17], this implies that $U^n$ is progressively measurable with respect to the filtration $\{\mathcal{F}^n_t\}$ (see Proposition 1.13 of [12]). We now assume that the primitive processes satisfy a functional strong law of large numbers. Recall that the abbreviation u.o.c. represents uniform convergence on compact time intervals.

ASSUMPTION 4.1.

1. There exists $\bar{u} \in \mathbb{R}^J_+$ such that a.s.

$$\lim_{n \to \infty} \frac{U^n(0)}{n} = \bar{u}.$$

2. For each $n \in \mathbb{N}$, there exists $\gamma^n \in \mathbb{R}^J_+$ such that a.s.

$$\lim_{m \to \infty} \frac{H^n(mt)}{m} = \gamma^n t,$$

where the convergence is u.o.c.

3. There exists $\gamma \in \mathbb{R}^J_+$ such that

$$\lim_{n \to \infty} \gamma^n = \gamma.$$

Recall that $\iota : \mathbb{R}_+ \to \mathbb{R}_+$ is the identity map. Define

$$\nu \doteq \gamma - \alpha, \tag{4.2}$$

and note that for $j \in \mathcal{I}$, $-\nu_j$ represents the amount of nominal service capacity allocated to class $j$ that is in excess of its mean arrival rate. Also, let

$$\overline{X} \doteq \bar{u} + \nu\iota, \qquad \overline{U} \doteq \Gamma(\overline{X}), \qquad \overline{Y} \doteq \overline{U} - \overline{X} \tag{4.3}$$

and

$$\overline{T} \doteq \alpha\iota - \overline{Y} = \gamma\iota + \bar{u} - \overline{U}. \tag{4.4}$$

The following result was established in [17].

THEOREM 4.2 (Fluid limits for the GPS model). *Suppose Assumptions 2.1 and 4.1 hold, and let $\overline{U}$ and $\overline{Y}$ be given by (4.3). Then $P$ a.s., as $n \to \infty$, $\overline{U}^n \to \overline{U}$, $\overline{Y}^n \to \overline{Y}$ and $\overline{T}^n \to \overline{T}$ u.o.c. Moreover, if $\sum_{i=1}^J \gamma_i = 1$ and $\bar{u} = 0$, then $\overline{U} = 0$ and $\overline{T} = \gamma \iota$.*

PROOF. The first statement is Theorem 4.3 of [17]. The condition $\sum_{i=1}^J \gamma_i = 1$ implies that $\sum_{i=1}^J \nu_i = 0$ and, thus, the second statement follows from Lemma 4.4(2) of [17]. □

REMARK 4.3. In Theorem 4.2 the fluid limit of the unfinished work is represented in terms of the GPS SP: $\overline{U} = \Gamma(\overline{X})$. As shown in [7], it is also possible to equivalently represent the fluid limit as the unique solution to a system of coupled differential equations. While the latter may in some sense provide a more intuitive description of the fluid limit, it does not extend to more general continuous inputs, as is required for the heavy traffic analysis. In contrast, the use of the GPS SP (ESP) provides a unified framework in which to study the pre-limit, fluid limit and diffusion limit. It is therefore natural and more convenient to work throughout with the SP formulation. Indeed, this level of abstraction allows one to better understand the connection between the nature of reallocation of service (as embodied in the directions of constraint) and the continuity and monotonicity properties of the map that are used in establishing the limit theorems. As a result, we expect that this approach may be more readily generalizable to other situations, including those in which a simple differential equation characterization of the fluid limit is not available.

4.2. *The invariant manifold of the fluid limit.* The goal of this section is to identify the so-called invariant manifold of the fluid limit. We first consider the task of identifying the set of strictly subcritical classes, namely, the sources whose long-run average arrival rate $\gamma_j$ is strictly less than the long-run average service rate available to them (after redistribution of service by the GPS discipline) for all sufficiently large $t$. For SPs associated with other queueing networks that have been studied in the literature (see, e.g., [4] or [21]), the linear independence of the associated constraint directions helps simplify this task. For example, in the study of open single-class queueing networks in [4], the associated unfinished work $U^*$ is represented as the image of the corresponding netput process $X^*$ under the associated SM $\Gamma^*$: $U^* = \Gamma^*(X^*) = X^* + Y^*$, where $Y^*$ now admits the decomposition $Y^* = D\xi^*$, where $D = (I - P')$ is a nonsingular matrix, $I$ is the identity matrix, $P$ is the so-called routing matrix of the network, the $i$th column of $D$ represents the direction of constraint associated with the face $\{x_i = 0\}$, and $\xi^*$ is the component-wise nondecreasing process that characterizes the cumulative idle time or excess capacity at each station (note that, somewhat unfortunately, the term $\xi^*$ here is denoted by the letter $Y$ in [4]). The fluid limit $\bar{\xi}^*$ of the vector of cumulative idleness admits the explicit expression $\bar{\xi}^* = D^{-1}\overline{Y}^*$, where $\overline{Y}^* = \overline{U}^* - \overline{X}^*$, with $\overline{U}^* = \Gamma^*(\overline{X}^*)$ being a continuous functional of $\overline{X}^*$ that is identically zero under the overall heavy traffic condition (the bar quantities here all refer to the fluid limits of the original quantities). In this case, the set of strictly subcritical classes is precisely the set of classes $j$ for which $\theta_j^* \doteq d\bar{\xi}_j^*/dt$ is strictly positive for all sufficiently large $t$. Since $\theta_j^*$ can be explicitly recovered from $\overline{X}^*$, which is itself known explicitly in terms of the primitives, this simplifies the determination of the strictly subcritical classes.

In contrast, as shown in (3.1), the GPS directions of constraint are linearly dependent and, thus, such a simple linear algebraic relation between the netput and the cumulative idle time processes no longer holds in general. Nevertheless, using the geometric properties established in Lemma 3.1, we show below in Lemma 4.4 that an analogous decomposition into component processes is possible for *fluid trajectories*. However, it is important to note that this decomposition does not hold for arbitrary trajectories (see Lemma 5.3), thus necessitating a more careful analysis of the diffusion limit (see Section 5.1).

LEMMA 4.4 (A representation lemma). *Given $\bar{u} \in \mathbb{R}_+^J$ and $\nu \in \mathbb{R}^J$, there exists $\varepsilon > 0$ and $\chi \in \mathbb{R}^J$ such that*

$$\Gamma(\bar{u} + \nu\iota)(t) = \bar{u} + \chi t \qquad \text{for } t \in [0, \varepsilon).$$

*Moreover, the following two properties are satisfied:*

(i) *There exists a unique vector $\theta \in \mathbb{R}_+^{J+1}$ that satisfies*

$$\chi = \nu + \sum_{k=1}^{J+1} \theta_k d_k \tag{4.5}$$

*and*

$$\mathcal{I} \neq \{k \in \mathcal{I} : \theta_k > 0\} \subseteq \{k \in \mathcal{I} : \bar{u}_k = \chi_k = 0\}. \tag{4.6}$$

*Moreover, $\theta_{J+1} > 0$ if and only if $\bar{u} = \chi = 0$ and $\sum_{j=1}^{J} \nu_j < 0$.*

(ii) *If $\bar{u} = 0$, then $\chi = \pi(\nu)$, $\varepsilon$ can be chosen to be $\infty$ and $[\sum_{j=1}^J \chi_j] = [\sum_{j=1}^J \nu_j] \vee 0$.*

(iii) *If $\sum_{j=1}^J \nu_j \geq 0$, then*

$$0 \leq \sum_{k:\theta_k>0} \frac{\theta_k}{1-\beta_k} = \frac{-\sum_{k:\theta_k>0} \nu_k}{1 - \sum_{k:\theta_k>0} \beta_k} \tag{4.7}$$

*and, for every* $j \in \mathcal{I}$,

$$(4.8) \qquad \nu_j = \begin{cases} -\dfrac{\theta_j}{1-\beta_j} + \beta_j \left( \dfrac{-\sum_{k:\theta_k>0} \nu_k}{1-\sum_{k:\theta_k>0} \beta_k} \right), & \text{if } \theta_j > 0, \\ \chi_j + \beta_j \left( \dfrac{-\sum_{k:\theta_k>0} \nu_k}{1-\sum_{k:\theta_k>0} \beta_k} \right), & \text{if } \theta_j = 0. \end{cases}$$

PROOF. Fix $\nu \in \mathbb{R}^J$ and let $\mathcal{E} \doteq \{k \in \mathcal{I} : \bar{u}_k = \chi_k = 0\}$. The existence of $\chi \in \mathbb{R}^J$ and $\varepsilon > 0$ results from the fact that images of affine trajectories under the GPS SM are piecewise affine (which was proved in Lemma 4.4 of [17]). Next, note that for all $t \in (0, \varepsilon)$, the definition of the GPS SP [in particular, the relation (2.4) and properties 4 and 5 of Definition 2.2] implies that $\chi$ is the unique vector that satisfies

$$(4.9) \qquad \chi - \nu = \frac{1}{t}[\Gamma(\bar{u} + \nu\iota)(t) - \bar{u} - \nu t] \in \begin{cases} \text{cone}[d_j, j \in \mathcal{E}], \\ \quad \text{if } \mathcal{E} \neq \mathcal{I}, \\ \text{cone}[d_j, j \in \mathcal{I}, d_{J+1}], \\ \quad \text{if } \mathcal{E} = \mathcal{I}. \end{cases}$$

This immediately ensures the existence of a vector $\theta \in \mathbb{R}_+^{J+1}$ that satisfies (4.5) and (4.6). To show uniqueness of $\theta$, first taking inner products of (4.5) with $d_{J+1}$ shows that

$$(4.10) \qquad \langle \chi, d_{J+1} \rangle = \langle \nu, d_{J+1} \rangle + \theta_{J+1}.$$

If $\mathcal{E} \neq \mathcal{I}$, then (4.9) and Lemma 3.1 show that $\chi - \nu \in H$, and so $\theta_{J+1} = 0$ by (4.10). On the other hand, if $\mathcal{E} = \mathcal{I}$, then (4.10) uniquely determines $\theta_{J+1}$, $\theta_{J+1} > 0$ if and only if $\sum_{j \in \mathcal{I}} \nu_j < 0$ and $\chi - \nu - \theta_{J+1} d_{J+1} \in H$. The last two statements, when combined with Lemma 3.1(ii), establish uniqueness of the representation (4.5) and the condition on $\theta_{J+1}$, thus proving property (i).

Now suppose that $\bar{u} = 0$. Then a simple consequence of the definition (3.6) of the GPS projection is that $\pi(\nu)\iota$ solves the SP for $\nu\iota$. Therefore, by uniqueness of solutions to the GPS SP, $\chi = \pi(\nu) \in \mathbb{R}_+^J$ and $\varepsilon$ can be chosen to be $\infty$. The second relation in property (ii) can be deduced in a straightforward manner from the fact that $\chi = \pi(\nu) \in \mathbb{R}_+^J$ and the last statement of property (i). The latter also shows that when $\sum_{j=1}^J \nu_j \geq 0$, $\theta_{J+1} = 0$ and $w \doteq \chi - \nu \in H$. When combined with (3.3) and (3.4) of Lemma 3.1, elementary algebra immediately yields property (iii). □

Recall from (4.2) that $\nu = \gamma - \alpha$. When $\sum_{j \in \mathcal{I}} \nu_j < 0$, it is intuitively clear that all classes will be strictly subcritical. The interesting case is thus when $\sum_{j \in \mathcal{I}} \nu_j \geq 0$. Suppose, in addition, that $\bar{u} = 0$, and note that by Lemma 4.4(ii), in this case $\chi = \pi(\nu) \in \mathbb{R}_+^J$. Let $\theta \in \mathbb{R}_+^{J+1}$ be the unique vector in the representation (4.5) for $\chi = \pi(\nu)$. With reference to the discussion prior to Lemma 4.4, it is natural to introduce the following definition:

$$(4.11) \qquad \mathcal{S}_0(\nu) \doteq \{j \in \mathcal{I} : \theta_j > 0\}.$$

REMARK 4.5. We now present an alternative characterization of $\mathcal{S}_0(\nu)$ that may appear more intuitive to some readers. When $\sum_{j \in \mathcal{I}} \nu_j \geq 0$, we claim that the set $\mathcal{S}_0(\nu)$ defined in (4.11) also admits the following alternative characterization: $\mathcal{S}_0(\nu)$ is the unique set $S$ that satisfies

$$(4.12) \qquad \sum_{j \in S} (\alpha_j - \gamma_j) > 0$$

and

$$(4.13) \quad \gamma_j < \alpha_j + \frac{\beta_j}{1 - \sum_{k \in S} \beta_k} \sum_{k \in S} (\alpha_k - \gamma_k) \quad \text{if and only if } j \in S.$$

Since the strictly subcritical classes receive no extra capacity from the other classes, a necessary condition for the set $S$ to be strictly subcritical is that the sum of the nominal capacities available to all classes in $S$ is strictly larger than the sum of the mean arrival rates of classes in $S$—this leads to the first condition (4.12). Moreover, it is precisely this excess that is redistributed to the remaining classes, with class $j$ receiving a fraction proportional to $\beta_j$, for $j \notin S$. This leads naturally to the "only if" part of the condition (4.13), which simply states that for any class that is not strictly subcritical, the total (reallocated) service capacity available to it is no greater than its mean (long-run) arrival rate. The "if" part of (4.13) is not as straightforward to justify a priori. Nevertheless, as shown below, it turns out to be the correct additional condition that uniquely characterizes $\mathcal{S}_0(\nu)$.

To see that (4.12) and (4.13) uniquely characterize $\mathcal{S}_0(\nu)$, let $S$ be any set satisfying (4.12) and (4.13). Since $\sum_{j \in \mathcal{I}} \nu_j = \sum_{j \in \mathcal{I}} (\gamma_j - \alpha_j) \geq 0$, the relation (4.12) ensures that $S \neq \mathcal{I}$ and thus the relation (4.13) is well defined. Now, the two conditions above are equivalent to the relations

$$(4.14) \qquad r_S \doteq \frac{-\sum_{j \in S} \nu_j}{1 - \sum_{j \in S} \beta_j} > 0$$

and

$$(4.15) \qquad \frac{\nu_j}{\beta_j} < r_S \qquad \text{iff } j \in S.$$

We now show that uniqueness of the set $S$ is essentially equivalent to uniqueness of the representation (4.5). Suppose that we are given a set $S \neq \mathcal{I}$ satisfying (4.14) and (4.15). Define

$$\tilde{\theta}_j \doteq -(1 - \beta_j)\nu_j + \beta_j(1 - \beta_j)r_S \qquad \text{for } j \in S.$$

Then elementary algebra shows that (4.14) and (4.15) imply that $\tilde{\theta}_j > 0$ for $j \in S$, and that $\tilde{\kappa} \doteq \nu + \sum_{j \in S} \tilde{\theta}_k d_k$ satisfies $\tilde{\kappa} = \pi(\nu)$. However, uniqueness of $\pi$ dictates that $\tilde{\kappa} = \kappa$. Uniqueness of the representation (4.5) then shows that $\tilde{\theta} = \theta$, with $\theta$ as in (4.11), and hence that $S = \mathcal{S}_0(\nu)$.

Now, let the sources be ordered so that

$$(4.16) \qquad \frac{\nu_1}{\beta_1} \geq \frac{\nu_2}{\beta_2} \geq \cdots \geq \frac{\nu_J}{\beta_J}.$$

It then follows immediately from Remark 4.5—in particular, relation (4.15)— that when $\sum_{j \in \mathcal{I}} \nu_j \geq 0$ either $\mathcal{S}_0(\nu) = \varnothing$ or $\mathcal{S}_0(\nu)$ has the form $\{j_*, \ldots, J\}$ for some $j_* \in \mathcal{I}$. Given $\nu \in \mathbb{R}^J$, we now define the set of strictly subcritical classes to be

$$(4.17) \qquad \mathcal{S}(\nu) = \begin{cases} \mathcal{I}, & \text{if } \sum_{j \in \mathcal{I}} \nu_j < 0, \\ \mathcal{S}_0(\nu), & \text{if } \sum_{j \in \mathcal{I}} \nu_j \geq 0, \end{cases}$$

and let $\mathcal{M}(\nu)$ be defined as follows:

$$(4.18) \qquad \mathcal{M}(\nu) \doteq \{x \in \mathbb{R}_+^J : x_i = 0 \text{ for every } i \in \mathcal{S}(\nu)\}.$$

Note that $\mathcal{M}(\nu) = \{0\}$ is 0-dimensional if $\sum_{j \in \mathcal{I}} \nu_j < 0$ and $\mathcal{M}(\nu) = \mathbb{R}_+^J$ is $J$-dimensional if $\sum_{j \in \mathcal{I}} \nu_j \geq 0$ and $\min_j \nu_j \geq 0$. The next result shows that $\mathcal{M}(\nu)$ acts as an invariant manifold for the fluid limit.

THEOREM 4.6 (Invariant manifold for the fluid limit). *Given $\nu \in \mathbb{R}^J$ and $\bar{u} \in \mathbb{R}_+^J$, let $\kappa \doteq \pi(\nu)$, $\overline{U} \doteq \Gamma(\bar{u} + \nu \iota)$, and let $\mathcal{S}(\nu)$ and $\mathcal{M}(\nu)$ be defined as in (4.17) and (4.18), respectively. Then the following properties are satisfied:*

(i) *If $\bar{u} \in \mathcal{M}(\nu)$, then $\overline{U}(t) = \bar{u} + \kappa t \in \mathcal{M}(\nu)$ for all $t \in (0, \infty)$.*

(ii) *Given any compact set $G \subset \mathbb{R}_+^J$, there exists $T = T(G) < \infty$ such that for every $\bar{u} \in G$, $\overline{U}(t) \in \mathcal{M}(\nu)$ for all $t \geq T$.*

(iii) *There exists $\bar{u} \in \mathcal{M}(\nu)$ such that $\overline{U}_i(t) > 0$ for every $t \in (0, \infty)$ and $i \in \mathcal{I} \setminus \mathcal{S}(\nu)$.*

(iv) *When $\sum_{j \in \mathcal{I}} \nu_j \leq 0$, $\mathcal{M}(\nu)$ admits the equivalent representation*

$$(4.19) \qquad \mathcal{M}(\nu) = \{\bar{u} \in \mathbb{R}_+^J : \Gamma(\bar{u} + \nu \iota)(t) = \bar{u} \text{ for all } t \in [0, \infty)\}.$$

(v) *Furthermore, if $\sum_{j \in \mathcal{I}} \nu_j = 0$ and $\min_{j \in \mathcal{I}} \nu_j < 0$, then there exists $j_* \in \mathcal{I}$ such that*

$$(4.20) \quad \frac{\nu_j}{\beta_j} = \frac{\nu_1}{\beta_1} \quad \text{for } j < j_* \quad \text{and} \quad \mathcal{S}(\nu) = \left\{j \in \mathcal{I} : \frac{\nu_j}{\beta_j} < \frac{\nu_1}{\beta_1}\right\}.$$

*Also,*

$$(4.21) \qquad \nu_j = \frac{-\beta_j (\sum_{j \in S(\nu)} \nu_j)}{1 - \sum_{k \in S(\nu)} \beta_k} \qquad \text{for } j \notin S(\nu).$$

PROOF. Fix $\nu \in \mathbb{R}^J$. When $\sum_{i \in \mathcal{I}} \nu_i < 0$, $\mathcal{S}(\nu) = \mathcal{I}$, $\mathcal{M}(\nu) = \{0\}$ and so the first two statements of the theorem follow directly from Lemma 4.4 of [17] (also see the results of [3]), and the third statement holds trivially. For the case when $\sum_{i \in \mathcal{I}} \nu_i \geq 0$ and $\min_{i \in \mathcal{I}} \nu_i \geq 0$, the first two statements of the theorem are a trivial consequence of the fact that $\mathcal{M}(\nu) = \mathbb{R}_+^J$ and the third statement is satisfied by any $\bar{u}$ with $\bar{u}_i > 0$ for all $i \in \mathcal{I}$, since for such $\bar{u}$, $\Gamma(\bar{u} + \nu\iota) = \bar{u} + \nu\iota$.

Therefore, for the rest of the proof of properties (i)–(iii), we shall assume that $\sum_{i \in \mathcal{I}} \nu_i \geq 0$ and $\min_{i \in \mathcal{I}} \nu_i < 0$, in which case $S(\nu) = S_0(\nu)$. We start by proving property (ii). Given any compact set $G \subset \mathbb{R}_+^J$, let $\bar{u}^* \in \mathbb{R}_+^J$ be such that for every $i \in \mathcal{I}$, $\bar{u}_i^* = \max_{\bar{u} \in G} \bar{u}_i$ and let $\overline{U}^* \doteq \Gamma(\bar{u}^* + \nu\iota)$. The comparison principle in Theorem 3.2 then guarantees that for any $\bar{u} \in G$, $\overline{U}_i(t) \leq \overline{U}_i^*(t)$ for every $i \in \mathcal{I}$, where $\overline{U} = \Gamma(\bar{u} + \nu\iota)$. Thus, in order to establish (ii), it suffices to show that there exists a finite $T < \infty$ such that $\overline{U}_j^*(t) = 0$ for every $j \in \mathcal{S}_0(\nu)$ and $t \geq T$. By the monotonicity property of the GPS SM established in Lemma 4.4 of [17], we know that

(4.22) $\{i \in \mathcal{I} : \overline{U}_i^*(t) = 0\} \subseteq \{i \in \mathcal{I} : \overline{U}_i^*(t+s) = 0\}$ for all $s, t > 0$.

As a result, there must exist $\mathcal{E} \subseteq \mathcal{I}$ and $T < \infty$ such that

(4.23) $\mathcal{E} = \{j \in \mathcal{I} : \overline{U}_j^*(t) = 0\}$ for all $t \geq T$.

Since $\overline{U}_i^*$ is piecewise affine, this implies that the slope of $\overline{U}^*(t)$ for all $t > T$ is well defined and equal to $\chi$, where $\chi$ satisfies $\chi_j = 0$ for $j \in \mathcal{E}$ and $\chi_j \geq 0$ for $j \notin \mathcal{E}$. On the other hand, by uniqueness of the SM, it is easy to see that, for $s \geq 0$,

(4.24) $\overline{U}^*(T+s) = \Gamma(\bar{u}^* + \nu\iota)(T+s) = \Gamma(\overline{U}^*(T) + \nu\iota)(s) = \overline{U}^*(T) + \chi s$.

From the definition of the SP, the last equality implies that for every $s \geq 0$,

$$(\chi - \nu) \in d(\overline{U}^*(T) + \chi s) = \text{cone}[d_i, i \in \mathcal{E}] \subseteq d(\chi).$$

Since $\chi \in \mathbb{R}_+^J$, by uniqueness of the projection $\pi$, this show that $\chi = \kappa \doteq \pi(\nu)$. The definition (4.11) of $\mathcal{S}_0(\nu)$ and relation (4.6) of Lemma 4.4 then shows that $\mathcal{S}_0(\nu) \subseteq \mathcal{E}$, which completes the proof of property (ii).

The definition of $\mathcal{S}_0(\nu)$ implies that $\kappa - \nu \in \text{cone}[d_i : i \in \mathcal{S}_0(\nu)] \subseteq \text{cone}[d_i : \kappa_i = 0]$. If $\bar{u} \in \mathcal{M}(\nu)$, then this implies that $\kappa - \nu \subseteq \text{cone}[d_i : \kappa_i = \bar{u}_i = 0]$, and thus, by (4.5) and (4.6) of Lemma 4.4, $\overline{U}(t) = \Gamma(\bar{u} + \nu\iota)(t) = \bar{u} + \kappa t$ for all sufficiently small $t$ (and, in fact, for all $t$ since $\kappa \in \mathbb{R}_+^J$). The above argument also shows that if $\bar{u} \in \mathcal{M}(\nu)$ with $\bar{u}_j > 0$ for every $j \notin \mathcal{S}_0(\nu)$, $\overline{U}_j(t) = \bar{u}_j + \kappa_j t > 0$ for all $t \geq 0$ and $j \notin \mathcal{S}_0(\nu)$. This establishes properties (i) and (iii).

Next, note that property (ii) shows that the set of invariant points of the fluid limit [described precisely by the right-hand side of (4.19)] must be contained in $\mathcal{M}$. If, in addition, $\sum_{j \in \mathcal{I}} \nu_j \leq 0$, then $\kappa = \pi(\nu) = 0$. Thus, property (i) shows that $\bar{u} \in \mathcal{M}$ implies $\overline{U}(t) = \Gamma(\bar{u} + \nu\iota)(t) = \bar{u}$ for every $t \in [0, \infty)$, which establishes the fourth property. In particular, this implies that $\chi = \kappa = 0$. Relation (4.20) is then an immediate consequence of relation (4.8) and the fact that $1 \notin \mathcal{S}_0(\nu)$ [since $\mathcal{S}_0(\nu) \neq \mathcal{I}$]. In turn, (4.20) implies that $\nu_j = \beta_j \nu_1 / \beta_1$ for $j \notin \mathcal{S}_0(\nu)$. Summing the last equality over $j \notin \mathcal{S}_0(\nu)$, and using the fact that $\sum_{k \in \mathcal{I}} \nu_k = 0$ and $\sum_{k \in \mathcal{I}} \beta_k = 1$, we obtain $\nu_1/\beta_1 = -\sum_{k \in \mathcal{S}_0(\nu)} \nu_k / (1 - \sum_{k \in \mathcal{S}_0(\nu)} \beta_k)$. When combined, the last two relations yield (4.21), thus completing the proof of the theorem. $\square$

REMARK 4.7. Although we do not use this later in the paper, it is also possible to provide a purely dynamical systems characterization of the invariant manifold $\mathcal{M}(\nu)$ defined in (4.18). Suppose, given $\nu \in \mathbb{R}^J$, an attractor $\mathcal{A}(\nu)$ for the fluid limit is defined to be any cone in $\mathbb{R}_+^J$ that satisfies (a) $\bar{u} \in \mathcal{A}(\nu)$ implies $\overline{U}(t) \doteq \Gamma(\bar{u} + \nu\iota)(t) \in \mathcal{A}(\nu)$ for every $t \in [0, \infty)$, and (b) $\bar{u} \notin \mathcal{A}(\nu)$ implies $\lim_{t \to \infty} \overline{U}(t) \in \mathcal{A}(\nu)$. Then it can be shown that the invariant manifold $\mathcal{M}(\nu)$ is the intersection of $\mathbb{R}_+^J$ with the affine hull of any attractor $\mathcal{A}(\nu)$ for the fluid limit. In the subcritical case it is easy to see from Theorem 4.6 that the invariant manifold $\mathcal{M}(\nu)$ is the unique attractor. Indeed, as shown in Theorem 4.6(iv), in this case the invariant manifold $\mathcal{M}(\nu)$ can equivalently be characterized as the collection of invariant points for the fluid limit, which is the "standard" definition of an invariant manifold given, for example, in [2] and [13]. However, while the standard definition is limited to the subcritical case, our definition is more general in that it also applies to the supercritical case, where it provides information on how trajectories escape to infinity. Indeed, in the supercritical case it is not hard to show that a set is an attractor if and only if it is a cone contained in $\mathcal{M}(\nu)$ that contains the ray $\{\kappa t, t \geq 0\}$ and has nonempty interior relative to $\mathcal{M}(\nu)$. Since we do not use this property later, we omit a rigorous proof of this statement.

**5. Diffusion approximations for the unbalanced GPS model.** For simplicity, we assume throughout this section that the classes are numbered so as to satisfy the ordering (4.16). As in the previous section, we consider a sequence of networks with associated processes $H^n, n \in \mathbb{N}$, that satisfy Assumptions 2.1 and 4.1. Recall the defining equations for the fluid limit processes $\overline{U}$, $\overline{X}$, $\overline{Y}$ and $\overline{T}$ given in (4.3) and (4.4), and consider the associated diffusion scaled processes defined by

(5.1)
$$\widehat{H}^n \doteq \sqrt{n}[\overline{H}^n - \gamma^n \iota], \qquad \widehat{U}^n \doteq \sqrt{n}[\overline{U}^n - \overline{U}],$$
$$\widehat{X}^n \doteq \sqrt{n}[\overline{X}^n - \overline{X}], \qquad \widehat{Y}^n \doteq \sqrt{n}[\overline{Y}^n - \overline{Y}],$$
$$\widehat{T}^n \doteq \sqrt{n}[\overline{T}^n - \overline{T}].$$

To prove the heavy traffic limit theorem for the unfinished work process, we first assume that the primitive sequence $\{H^n\}$ satisfies, in addition, a functional central limit theorem. Let $\gamma \in \mathbb{R}_+^J$ be the vector in Assumption 4.1(3) and let $\nu \doteq \gamma - \alpha$.

ASSUMPTION 5.1.

1. There exists a random variable $\hat{u} \in \mathcal{M}(\nu)$ such that $P$ a.s.,
$$\lim_{n \to \infty} \frac{U^n(0)}{\sqrt{n}} = \hat{u}.$$

2. As $n \to \infty$,
$$\widehat{H}^n \Rightarrow B,$$
where $B$ is a driftless $J$-dimensional Brownian motion with covariance matrix $M^H$.

3. There exists $\hat{c} \in \mathbb{R}^J$ such that
$$\lim_{n \to \infty} \sqrt{n}(\gamma^n - \gamma) = \hat{c}.$$

We will also assume the heavy traffic condition
$$\sum_{i=1}^{J} \gamma_i = 1. \tag{5.2}$$

REMARK 5.2. (a) Note that Assumption 5.1(3) implies, in particular, that $\gamma^n \to \gamma$ as $n \to \infty$.

(b) By the Skorokhod representation theorem (see, e.g., Theorem 1.8 of [10]), there exists a probability space $(\tilde{\Omega}, \tilde{F}, \tilde{P})$, on which are defined $\mathcal{D}([0, \infty) : \mathbb{R}^J)$-valued random variables $\tilde{H}^n$, $n \in \mathbb{N}$, and $\tilde{B}$ such that $\tilde{H}^n \stackrel{d}{=} \widehat{H}^n$, $\tilde{B} \stackrel{d}{=} B$ and $\tilde{H}^n \to B$ a.s. uniformly on compact sets (u.o.c.). By an abuse of notation, we simply take Assumption 5.2(2) to mean that $\tilde{H}^n \to B$ a.s. u.o.c. Of course, what we ultimately prove is weak convergence and not a.s. convergence, and this is reflected in the theorem statements.

In the next section we summarize our approach to the GPS diffusion limit, and discuss its connection with related work. This section can be safely skipped without loss of continuity.

5.1. *General approach to the unbalanced GPS diffusion limit.* The fluid limit result summarized in Theorem 4.2 shows that under the heavy traffic condition (5.2), $\overline{U} = 0$. By (5.1), (2.5) and basic homogeneity properties of the SP, it then follows that
$$\widehat{U}^n = \sqrt{n}\overline{U}^n = \Gamma(\sqrt{n}\overline{X}^n) = \overline{\Gamma}(\sqrt{n}\overline{X}^n),$$

where $\overline{\Gamma}$ is the GPS ESM and the last equality follows by Remark 2.4. If $\sqrt{n}\overline{X}^n$ could be shown to converge to a limit $\widehat{X}$, then, since the GPS ESM is Lipschitz continuous, the continuous mapping theorem would immediately yield convergence of $\widehat{U}^n$ to $\overline{\Gamma}(\widehat{X})$. Combining (2.7) with Assumption 5.1, it is not hard to see that $\sqrt{n}\overline{X}^n$ converges if and only if $\gamma = \alpha$. This explains why the standard continuous mapping approach works only in the balanced case (see [17]). Indeed, in the (truly) unbalanced case, $\sqrt{n}\overline{X}^n$ certainly diverges because the long-run average arrival rate $\gamma_j$ for at least one critical class must be strictly greater than its nominal capacity. This class becomes critical (in the fluid limit) only because it receives extra service capacity from the strictly subcritical classes. Thus, the methods of [17], which carried out an analysis of the balanced GPS model, are not sufficient to analyze the more general, unbalanced case considered here.

A similar generalization was considered in the context of single-class queueing networks. Specifically, the diffusion analysis in [19] [where each station was assumed to be in heavy traffic—see condition (24) therein] was extended in [4] to the unbalanced case. As mentioned in Section 4.2, the unfinished work $U^*$ for the model in [4] is represented, as in this paper, in the form $U^* = X^* + Y^*$, where $X^*$ is the associated netput process, but where $Y^*$ now admits the decomposition $Y^* = D\xi^*$ for a nonsingular matrix $D$ and $\xi^*$ is the component-wise nondecreasing process that characterizes the cumulative idle time or excess capacity at each station. The $i$th column of $D$ in [4] is analogous to the direction of constraint $d_i$ in this paper, and the mapping from $X^*$ to $U^*$ corresponds to the SM $\Gamma^*$ considered here. The diffusion analysis in [4] strongly uses (i) the explicit representation $Y^* = D\xi^* = \sum_{i=1}^{J} \xi_i^* d_i$ (see the block decompositions in the statement of Theorem 6.1 of [4]), as well as (ii) the continuity of the mapping that takes $X^*$ to $\xi^*$, which is referred to there as the regulator mapping (see (3.3D)–(3.3E) and the discussions following (4.25) and (4.29) in [4]).

On the other hand, no such explicit representation is available for the GPS SP—in fact, the directions of constraint are linearly dependent, as shown in (3.1). Nevertheless, as shown in Lemma 4.4, an analogous decomposition holds when the GPS SM acts on affine trajectories. Indeed, suppose $\psi$ is an affine trajectory, say, of the form $\bar{u} + \nu\iota$ for some $\bar{u} \in \mathbb{R}_+^J$ and $\nu \in \mathbb{R}^J$, and $\phi = \Gamma(\psi)$, where $\Gamma$ is the GPS SM. Since $\phi$ is piecewise affine (i.e., with a finite number of changes of slope), from Lemma 4.4 and the property (4.24) of the SP one can infer that there exist unique measurable functions $\theta : [0,\infty) \to \mathbb{R}_+^{J+1}$ with the property that for a.e. $s \in [0,\infty)$,

$$\dot{\phi}(s) = \dot{\psi}(s) + \sum_{j=1}^{J+1} \theta_j(s) d_j,$$

where, for $j \in \mathcal{I}$,

(5.3) $\qquad \{j : \theta_j(s) > 0\} \subseteq \{j : \phi_j(s) = 0\}$

and $\theta_{J+1}(s) > 0$ implies $\phi(s) = 0$. Another application of the property (4.24) shows in fact that this also holds for piecewise affine $\psi$. Let $\xi(t) = \int_0^t \theta(s) \, ds$ for $t \in [0, \infty)$, and let $\Theta$ be the mapping defined on piecewise affine, continuous functions that takes $\psi$ to the corresponding $\xi$. In view of the methods used in [4] discussed above, it is natural to ask whether $\Theta$ can be extended to define a continuous mapping on the space of all continuous functions. The following lemma answers this question in the negative—in fact, it shows that there is no canonical way to extend the definition of $\Theta$ to all continuous trajectories. This constitutes a structural difference between the GPS SP and other SPs considered in the literature in the context of queueing networks.

LEMMA 5.3. *The mapping $\Theta$ cannot be extended to define a continuous mapping on the space of all continuous functions.*

PROOF. Let $\Gamma$ be the 2-dimensional GPS SM, with directions of constraint $d_1 = (1, -1), d_2 = (-1, 1)$ and $d_3 = (1, 1)/\sqrt{2}$. If there exists a continuous extension of $\Theta$ to all continuous trajectories, then by the Cauchy property, for any pair of sequences of piecewise affine continuous functions $\{\psi^{(1,n)}\}_{n \in \mathbb{N}}$ and $\{\psi^{(2,n)}\}_{n \in \mathbb{N}}$ such that $\psi^{(1,n)} - \psi^{(2,n)} \to 0$, we must have $\Theta(\psi^{(1,n)}) - \Theta(\psi^{(2,n)}) \to 0$. We show that this is not the case by constructing a counterexample. Define $\psi^{(1,n)} \equiv 0$, let $\psi^{(2,n)}(0) = 0$, and let

$$\dot{\psi}^{(2,n)} \doteq \begin{cases} \left(\dfrac{2^n}{n} d_1\right), & \text{if } t \in [2(m-1)2^{-n}, (2m-1)2^{-n}), \\ \left(\dfrac{2^n}{n} d_2\right), & \text{if } t \in [(2m-1)2^{-n}, (2m)2^{-n}) \end{cases}$$

for $m = 1, \ldots, 2^{n-1}$. Then it is easy to see that

$$\sup_{s \in [0,1]} |\psi^{(1,n)}(s) - \psi^{(2,n)}(s)| = \sup_{s \in [0,1]} |\psi^{(2,n)}(s)| \leq \sqrt{2}/n,$$

and thus, $\psi^{(1,n)} - \psi^{(2,n)} \to 0$. On the other hand, if $\xi^{(1,n)} \doteq \Theta(\psi^{(1,n)})$ and $\xi^{(2,n)} \doteq \Theta(\psi^{(2,n)})$, it is easily verified that $\xi^{(1,n)} \equiv 0$, while $\xi_i^{(2,n)}(1) = 2^{n-1}/n$ for $i = 1, 2$ and $\xi_3^{(2,n)}(1) = 0$. This implies that

$$\lim_{n \to \infty} |\xi^{(2,n)}(1) - \xi^{(1,n)}(1)| \to \infty,$$

which completes the proof. $\square$

The aim of the above discussion was to explain why the analysis of the unbalanced GPS model does not fall into any of the previously existing frameworks for establishing diffusion approximations of queueing networks. We now briefly describe the approach taken in this paper to establishing diffusion limits, which entails first showing that there is state-space collapse, in the sense that the strictly sub-critical classes vanish in the diffusion limit, and then using this information to provide a nice characterization of the behavior of the remaining, critical classes. Lemma 5.3 suggests that it would be preferable to work, as far as possible, in the pre-limit since the GPS mapping is better behaved on the pre-limit than on the limit functions. Thus, we first introduce a modified work arrival process in Section 5.2 and identify the diffusion limit of the unfinished work associated with this work arrival process. The modified work arrival process is second-order equivalent to the original one in the sense that their diffusion limits coincide (see Theorem 5.4). However, as shown in Lemma 5.5, the advantage of working with the modified arrival process is that the state-space collapse takes place in the pre-limit itself. This facilitates a simple characterization of the critical classes in terms of a certain reduced map, which is introduced in Section 5.3. The reduced map can be characterized as a GPS SM on a lower-dimensional space (associated with the critical classes) with appropriately modified weights and is thus continuous. These results are finally combined with the comparison principle of Section 3.2 in order to establish state-space collapse and obtain an explicit characterization of the diffusion limit in Section 5.4. The diffusion limit identifies precisely how the covariance structure of the unfinished work of the critical classes is influenced by the variance of the cumulative work arrival processes of the strictly subcritical classes.

5.2. *A modified work arrival process.* In this section we introduce a sequence of modified cumulative work arrival processes $\{\mathcal{H}^n, n \in \mathbb{N}\}$ obtained by smoothing the sequence of class $j$ arrival process $\{H_j^n, n \in \mathbb{N}\}$ for $j \in \mathcal{S}(\nu)$. In Theorem 5.4 below, we show that $\{\mathcal{H}^n, n \in \mathbb{N}\}$ is equal to $\{H^n, n \in \mathbb{N}\}$ up to second order in the sense that their diffusion limits coincide. This is convenient because, as shown in the next section, the limit of the sequence of unfinished work processes associated with the modified arrival processes can be obtained by applying Theorem 5.6.

A rigorous definition of the sequence $\{\mathcal{H}^n, n \in \mathbb{N}\}$ is given below. Let $\{\varepsilon^n, n \in \mathbb{N}\}$ be a sequence of positive numbers such that

$$(5.4) \qquad \lim_{n \to \infty} \varepsilon^n = 0 \quad \text{and} \quad \lim_{n \to \infty} \sqrt{n}\varepsilon^n = \infty.$$

For $n \in \mathbb{N}$, define

$$(5.5) \qquad \tilde{\gamma}_j^n \doteq \begin{cases} \gamma_j^n + \varepsilon^n, & \text{for } j \in \mathcal{S}(\nu), \\ \gamma_j^n, & \text{for } j \in \mathcal{I} \setminus \mathcal{S}(\nu) \end{cases}$$

and $\tilde{\nu}^n \doteq \tilde{\gamma}^n - \alpha$. For $j \notin \mathcal{S}(\nu)$, let $\mathcal{H}_j^n = H_j^n$ and for $j \in \mathcal{S}(\nu)$, let $\mathcal{H}_j^n$ be the departure process from a queue that is initially empty, has cumulative work arrival process $H_j^n$ and a deterministic service rate of $\tilde{\gamma}_j^n$. In other words, we can write

$$(5.6) \qquad \mathcal{H}_j^n = \begin{cases} H_j^n - \Gamma_1(H_j^n - \tilde{\gamma}_j^n \iota), & \text{for } j \in \mathcal{S}(\nu), \\ H_j^n, & \text{for } j \notin \mathcal{S}(\nu), \end{cases}$$

where $\Gamma_1$ is the one-dimensional SM defined in (3.9). It should be noted that $\mathcal{H}_j^n$ does not satisfy Assumption 2.1(2). Nonetheless, it plays an essential role in proving our heavy traffic limit theorem. Let the corresponding fluid scaled and diffusion scaled processes $\overline{\mathcal{H}}^n$ and $\widehat{\mathcal{H}}^n$ be defined as in (4.1) and (5.1), respectively, with $H$ replaced by $\mathcal{H}$. The following theorem states a second-order equivalence result between the sequences $\{H^n\}$ and $\{\mathcal{H}^n\}$.

THEOREM 5.4 (Properties of the modified arrival sequence). *Suppose Assumptions 4.1 and 5.1 are satisfied. Then the sequence $\{\mathcal{H}^n\}$ satisfies the following properties:*

1. $\mathcal{H}_j^n(t) - \mathcal{H}_j^n(s) \leq \tilde{\gamma}_j^n(t-s)$ for any $0 \leq s \leq t$ and every $j \in \mathcal{S}(\nu)$.
2. $\overline{\mathcal{H}^n}^m \to \gamma^n \iota$ u.o.c. as $m \to \infty$, where $\overline{\mathcal{H}^n}^m(t) \doteq \frac{\mathcal{H}^n(mt)}{m}$.
3. $\widehat{\mathcal{H}}^n \Rightarrow B$ as $n \to \infty$.

PROOF. From (5.6) and the explicit expression (3.9) for the one-dimensional SM we see that for $t \in [0, \infty)$ and $j \in \mathcal{S}(\nu)$,

$$(5.7) \qquad \mathcal{H}_j^n(t) = \tilde{\gamma}_j^n t - \sup_{s \in [0,t]} [\tilde{\gamma}_j^n s - H_j^n(s)] \vee 0,$$

from which the first property immediately follows. Next, note that by (5.6), the definitions (5.5) of $\tilde{\gamma}$ and (5.1) of $\widehat{H}^n$, and elementary scaling properties of the one-dimensional SM, we have

$$\overline{\mathcal{H}_j^n}^m = \overline{H_j^n}^m - \Gamma_1(\overline{H_j^n}^m - \gamma_j^n \iota - \varepsilon^n \iota) \quad \text{and} \quad \widehat{\mathcal{H}}_j^n = \widehat{H}_j^n - \Gamma_1(\widehat{H}_j^n - \sqrt{n}\varepsilon^n \iota).$$

Taking limits as $m \to \infty$ in the first equation above and using the first limit in (5.4), along with Assumption 4.1(2), the fact that $\Gamma_1$ is continuous and $\Gamma_1(0) = 0$, we obtain the second property of the theorem.

By Peterson's "crushing lemma" ([15], Lemma 2), $\Gamma_1(\widehat{H}_j^n - \sqrt{n}\varepsilon^n \iota) \to 0$. [Although Peterson assumes that $\varepsilon^n \geq \epsilon > 0$, his proof goes through with $\varepsilon^n \to 0$ if $\sqrt{n}\varepsilon^n \to \infty$. Also, Peterson assumes mutually independent i.i.d. sequences of interarrival and service times, but his proof goes through under our Assumption 5.1(2).] The third property is then an immediate consequence of Assumption 5.1(2) and Remark 5.2. $\square$

For $n \in \mathbb{N}$, define

$$\mathcal{U}_j^n(0) = \begin{cases} 0, & \text{if } j \in \mathcal{S}(\nu), \\ U_j^n(0), & \text{if } j \notin \mathcal{S}(\nu) \end{cases}$$

and let $\mathcal{X}^n \doteq \mathcal{U}^n(0) + \mathcal{H}^n - \alpha \iota$. Moreover, let $\mathcal{U}^n \doteq \Gamma(\mathcal{X}^n)$.

The next lemma identifies the property that makes the unfinished work associated with the modified arrival sequence $\{\mathcal{H}^n, n \in \mathbb{N}\}$ easier to analyze than that associated with the original arrival sequence $\{H^n, n \in \mathbb{N}\}$. Indeed, it establishes state-space collapse for the pre-limit sequence of unfinished work processes associated with the sequence of modified arrival processes.

LEMMA 5.5. *Suppose Assumptions* 4.1, 5.1 *and the overall heavy traffic condition* (5.2) *are satisfied. Then for every* $\omega \in \Omega$, *there exists* $N = N(\omega) < \infty$ *such that, for every* $n \geq N(\omega)$, $\mathcal{U}^n(t)(\omega) \in \mathcal{M}(\nu)$ *for all* $t \in [0, \infty)$.

PROOF. Fix $\omega \in \Omega$ and drop the explicit dependence on $\omega$. Recall that $\tilde{\nu}^n = \tilde{\gamma}^n - \alpha$, and define $\tilde{\Phi}^n \doteq \Gamma(\mathcal{U}^n(0) + \tilde{\nu}^n \iota)$ and $\Phi^n \doteq \Gamma(\mathcal{U}^n(0) + \nu \iota)$. Then Theorem 5.4(1) along with the comparison principle of Theorem 3.2 shows that $\mathcal{U}^n(t) \leq \tilde{\Phi}^n(t)$ for $t \in [0, \infty)$. Thus, to prove the lemma, it suffices to demonstrate the existence of $N < \infty$ such that, for $n \geq N$, $\tilde{\Phi}^n(t) \in \mathcal{M}(\nu)$ for every $t \in [0, \infty)$. By Lemma 4.4, there exist $\delta^n > 0$ and $\chi^n, \tilde{\chi}^n \in \mathbb{R}^J$ such that, for $t \in [0, \delta^n)$,

$$\tilde{\Phi}^n(t) = \mathcal{U}^n(0) + \tilde{\chi}^n t \quad \text{and} \quad \Phi^n(t) = \mathcal{U}^n(0) + \chi^n t.$$

In addition, for every $n \in \mathbb{N}$, since $\mathcal{U}^n(0) \in \mathcal{M}(\nu)$, Theorem 4.6 shows that $\chi^n = \kappa = \pi(\nu)$. Moreover, for $n \in \mathbb{N}$, let $\tilde{\theta}^n \in \mathbb{R}_+^{J+1}$ be the unique vector in the representation (4.5) for $\tilde{\chi}^n$ and let $\theta \in \mathbb{R}_+^{J+1}$ be the corresponding unique vector in the representation (4.5) for $\kappa$. Lipschitz continuity of the GPS SM $\Gamma$ and the convergence $\tilde{\nu}_n \to \nu$ for $n \to \infty$, which is guaranteed by Assumption 4.1(3) and the first relation in (5.4), imply $\tilde{\chi}^n \to \kappa$ as $n \to \infty$. The uniqueness of the representation (4.5) then dictates that $\tilde{\theta}_n \to \theta$ as $n \to \infty$. Thus, there exists $N < \infty$ such that, for every $n \geq N$,

$$\mathcal{S}_0(\nu) = \{j \in \mathcal{I} : \theta_j > 0\} \subseteq \{j \in \mathcal{I} : \tilde{\theta}_j^n > 0\} \subseteq \{j \in \mathcal{I} : \mathcal{U}_j^n(0) = \tilde{\chi}_j^n = 0\},$$

where the first equality follows from definition (4.11) and the last inclusion is a result of relation (4.6) of Lemma 4.4. This implies that $\tilde{\chi}_j^n = 0$ for $j \in \mathcal{S}_0(\nu) = \mathcal{S}(\nu)$, so that $\tilde{\Phi}_j^n(t) \in \mathcal{M}(\nu)$ for $t \in [0, \tilde{\delta}^n)$. By the monotonicity property of the GPS SM, which was proved in Lemma 4.4 of [17] and is spelled out in (4.22) here, it then follows that $\tilde{\Phi}_j^n(t) \in \mathcal{M}(\nu)$ for every $t \in [0, \infty)$, which completes the proof of the theorem.  □

5.3. *A reduced representation for the SM.* We assume throughout this section that $\mathcal{S}$ is a subset of $\mathcal{I}$, $\mathcal{S} \neq \mathcal{I}$. The application of the main result, Theorem 5.6, of this section will be to the case when $\mathcal{S} = \mathcal{S}(\nu)$ is the set of strictly subcritical classes associated with some $\nu$ such that $\sum_{j=1}^J \nu_j = 0$. However, the results of this section are valid for an arbitrary set $\mathcal{S} \subset \mathcal{I}$. In order to formulate our reduced representation for the SM, we need the following notation. Let $\mathcal{K} \doteq \mathcal{I} \setminus \mathcal{S}$ and $K \doteq |\mathcal{K}| \geq 1$. Given a weight vector $\beta \in (0,1)^J$ (with $\sum_{i=1}^J \beta_i = 1$), let $\beta^\mathcal{K} = (\beta_k^\mathcal{K}, k \in \mathcal{K})$ be the $K$-dimensional vector defined by

$$(5.8) \qquad \beta_i^\mathcal{K} \doteq \frac{\beta_i}{\sum_{k \in \mathcal{K}} \beta_k} \qquad \text{for } i \in \mathcal{K}.$$

Also, for $i \in \mathcal{K}$, let $e_i^{\mathcal{K}}$ be the $K$-dimensional unit vector associated with the $i$th coordinate, which has 1 at the $i$th coordinate and 0 for all other coordinates in $\mathcal{K}$. Since $\mathcal{K} \neq \varnothing$ and $\beta \in (0,1)^J$, it follows that $\sum_{k \in \mathcal{K}} \beta_k > 0$. Hence, $\beta^{\mathcal{K}}$ is well defined, $\beta^{\mathcal{K}} \in (0,1)^K$ and $\sum_{k \in \mathcal{K}} \beta_k^{\mathcal{K}} = 1$, which shows that $\beta^{\mathcal{K}}$ is a $K$- dimensional weight vector. We now define the GPS SP associated with the weight vector $\beta^{\mathcal{K}}$. First, let the directions of constraint $\{d_i^{\mathcal{K}}, i \in \mathcal{K}, d_{K+1}^{\mathcal{K}}\}$ be defined in terms of the $K$-dimensional weight vector $\beta^{\mathcal{K}}$ in the same way that the directions of constraint $\{d_i, i \in \mathcal{I}, d_{J+1}\}$ were defined in terms of the $J$-dimensional weight vector $\beta$ in (2.3): more precisely, let $d_{K+1}^{\mathcal{K}} \doteq \sum_{i \in \mathcal{K}} e_i^{\mathcal{K}}/\sqrt{K}$ and

$$(5.9) \qquad d_i^{\mathcal{K}} \doteq e_i^{\mathcal{K}} - \sum_{k \in \mathcal{K} \setminus \{i\}} \frac{\beta_k^{\mathcal{K}} e_k^{\mathcal{K}}}{1 - \beta_i^{\mathcal{K}}} \qquad \text{for } i \in \mathcal{K}.$$

Recalling definition (5.8), for $i \in \mathcal{K}$, we can then write

$$(5.10) \qquad (d_i^{\mathcal{K}})_j \doteq \begin{cases} 1, & \text{if } j = i, \\ \dfrac{-\beta_j}{\sum_{k \in \mathcal{K}} \beta_k - \beta_i}, & \text{if } j \in \mathcal{K} \setminus \{i\}. \end{cases}$$

Define

$$\Lambda^{\mathcal{K}} \doteq \{(x_i, i \in \mathcal{K}): x_i \in \mathbb{R} \text{ for } i \in \mathcal{K}\},$$

$$\Lambda_+^{\mathcal{K}} \doteq \{(x_i, i \in \mathcal{K}): x_i \geq 0 \text{ for } i \in \mathcal{K}\}$$

and, analogously to the definition given in (2.4), let

$$d^{\mathcal{K}}(x) \doteq \left\{ \sum_{i \in I^{\mathcal{K}}(x)} a_i d_i^{\mathcal{K}} : a_i \geq 0 \text{ for } i \in I^{\mathcal{K}}(x) \right\},$$

where, for $x \in \Lambda_+^{\mathcal{K}}$, $I^{\mathcal{K}}(x)$ is defined as in (2.5), but with $\mathcal{I}$ replaced by $\mathcal{K}$ and $J$ by $K$. Similarly, let the "$\mathcal{K}$-reduced" GPS SP be as described in Definition 2.2, but with $\mathbb{R}^J$, $\mathbb{R}_+^J$ and $d(\cdot)$ replaced by $\Lambda^{\mathcal{K}}$, $\Lambda_+^{\mathcal{K}}$ and $d^{\mathcal{K}}(\cdot)$, respectively. Finally, let $\Gamma^{\mathcal{K}}$ be the associated SM, and note that it is well defined for the same reasons given in Remark 2.4. We now establish a correspondence between the original GPS SM and the $\mathcal{K}$-reduced SM.

THEOREM 5.6. *Given a weight vector $\beta \in (0,1)^J$, let $\Gamma$ be the associated GPS ESM, let $\psi \in \mathcal{D}([0,\infty); \mathbb{R}^J)$, let $\varphi = \overline{\Gamma}(\psi)$ and for $t \in [0,\infty)$, let $E(t) \doteq \{i \in \mathcal{I}: \varphi_i(t) = 0\}$. Moreover, suppose there exist $0 \leq T_0 < T < \infty$ and $\mathcal{S} \subsetneq \mathcal{I}$ such that $\mathcal{S} \subseteq E(t)$ for every $t \in [T_0, T)$. Then, if $\mathcal{K} \doteq \mathcal{I} \setminus \mathcal{S}$, the vector-valued process $\varphi^{\mathcal{K}} \doteq \{\varphi_i, i \in \mathcal{K}\}$ satisfies*

$$\varphi^{\mathcal{K}}(t) = \overline{\Gamma}^{\mathcal{K}}(\psi^{\mathcal{K}})(t) \qquad \text{for } t \in [T_0, T),$$

where $\overline{\Gamma}^{\mathcal{K}}$ is the $\mathcal{K}$-reduced GPS ESM defined above and $\psi^{\mathcal{K}} \doteq \{\psi_i^{\mathcal{K}}, i \in \mathcal{K}\}$ is given by

$$(5.11) \qquad \psi_i^{\mathcal{K}} \doteq \psi_i + \beta_i^{\mathcal{K}} \left( \sum_{j \in \mathcal{S}} \psi_j \right) \qquad \text{for } i \in \mathcal{K}.$$

Finally, if $\psi \in \text{dom}(\Gamma)$, then $\varphi^{\mathcal{K}} = \Gamma^{\mathcal{K}}(\psi^{\mathcal{K}})$.

PROOF. Since $\mathcal{K} \neq \varnothing$, $\beta^{\mathcal{K}}$ and $\overline{\Gamma}^{\mathcal{K}}$ are well defined. The theorem follows trivially if $\mathcal{K} = \mathcal{I}$. Thus, throughout the proof, we assume that $\mathcal{K} \subset \mathcal{I}$ or, equivalently, that $\mathcal{S} \neq \varnothing$. Since $\overline{\Gamma}$ is single-valued on its domain (see Remark 2.4), if $\varphi = \overline{\Gamma}(\psi)$, then $\tilde{\varphi} = \overline{\Gamma}(\tilde{\psi})$ for $t \in [0, \infty)$, where

$$\tilde{\varphi}(t) \doteq \varphi(T_0 + t) \quad \text{and} \quad \tilde{\psi}(t) \doteq \psi(T_0 + t) - \psi(T_0) + \varphi(T_0).$$

Thus, we can also assume without loss of generality that $T_0 = 0$.

For every $t \in [0, T)$, $\varphi(t) \in \mathbb{R}_+^J$ and, thus, $\varphi^{\mathcal{K}}(t) \in \Lambda_+^{\mathcal{K}}$. Define $\eta \doteq \varphi - \psi$ and $\eta^{\mathcal{K}} \doteq \varphi^{\mathcal{K}} - \psi^{\mathcal{K}}$. It remains to show that in the interval $[0, T)$, $(\varphi^{\mathcal{K}}, \eta^{\mathcal{K}})$ satisfy properties 3 and 4 of Definition 2.3 for the GPS ESP on $\Lambda_+^{\mathcal{K}}$ (corresponding to the weights $\beta^{\mathcal{K}}$). First, note that for $i \in \mathcal{K}$ and $t \in [0, T)$,

$$(5.12) \qquad \begin{aligned} \eta_i^{\mathcal{K}}(t) &= \eta_i(t) - \beta_i^{\mathcal{K}} \left( \sum_{k \in \mathcal{S}} \psi_k(t) \right) \\ &= \eta_i(t) + \beta_i^{\mathcal{K}} \left( \sum_{k \in \mathcal{S}} \eta_k(t) \right), \end{aligned}$$

where the second equality follows from the fact that $\varphi_k(t) = 0$ for $k \in \mathcal{S}$ and $\varphi = \psi + \eta$.

Now, fix $0 \leq s < t < T$. In analogy with the hyperplane $H$ defined in Lemma 3.1, let $H^{\mathcal{K}} \doteq \{x \in \Lambda^{\mathcal{K}} : \sum_{i \in \mathcal{K}} x_i = 0\}$. Recall that for every $x \in \mathbb{R}_+^J \setminus \{0\}$, $d(x) \subset H \subset d(0) = \{x \in \mathbb{R}^J : \sum_{i=1}^J x_i \geq 0\}$ and, likewise, for every $x \in \Lambda_+^{\mathcal{K}} \setminus \{0^{\mathcal{K}}\}$, $d^{\mathcal{K}}(x) \subset H^{\mathcal{K}} \subset d^{\mathcal{K}}(0^{\mathcal{K}}) = \{x \in \Lambda^{\mathcal{K}} : \sum_{i \in \mathcal{K}} x_i \geq 0\}$. Also, note that for any $u \in [0, T]$, $\sum_{i \in \mathcal{K}} \varphi_i^{\mathcal{K}}(u) = \sum_{i \in \mathcal{I}} \varphi_i(u)$ since $\varphi_i(u) = 0$ for $i \in \mathcal{S}$, and also that (5.11) implies $\sum_{i \in \mathcal{K}} \psi_i^{\mathcal{K}}(u) = \sum_{i \in \mathcal{I}} \psi_i(u)$. The last two statements along with the fact that $(\varphi, \eta)$ satisfy property 3 shows that

$$(5.13) \qquad \sum_{i \in \mathcal{K}} (\eta_i^{\mathcal{K}}(t) - \eta_i^{\mathcal{K}}(s)) = \sum_{i \in \mathcal{I}} (\eta_i(t) - \eta_i(s)) \geq 0,$$

with the inequality being replaced by an equality if, for every $u \in (s, t]$, $\varphi(u) \neq 0$. This shows, in particular, that $\eta^{\mathcal{K}}(t) - \eta^{\mathcal{K}}(s) \in d^{\mathcal{K}}(0^{\mathcal{K}})$ and $\eta^{\mathcal{K}}(t) - \eta^{\mathcal{K}}(s) \in H^{\mathcal{K}}$ if, for every $u \in (s, t]$, $\varphi(u) \neq 0$. Thus, this proves property 3 in the case when either $\varphi(u) = 0$ for some $u \in (s, t]$ or when $\varphi(u) \neq 0$ for every $u \in (s, t]$, but $\overline{\text{co}}[\bigcup_{u \in (s,t]} d^{\mathcal{K}}(\varphi^{\mathcal{K}}(u))] = H^{\mathcal{K}}$.

We now consider the remaining case when $\varphi(u) \neq 0$ for every $u \in (s, t]$, and $\overline{\operatorname{co}}[\bigcup_{u \in (s,t]} d^{\mathcal{K}}(\varphi^{\mathcal{K}}(u))] \subsetneq H^{\mathcal{K}}$. By Lemma 3.1, the latter relation implies that there exists $j \in \mathcal{K}$ such that, for every $u \in (s, t]$, $\varphi_j(u) = \varphi_j^{\mathcal{K}}(u) > 0$. Since $(\varphi, \eta)$ satisfy property 3 of the ESM and the vectors $d_i, i \in \overline{E} \doteq \mathcal{I} \setminus \{j\}$, are linearly independent, there exist unique $\theta_k \geq 0$, $k \in \overline{E}$, such that

$$\eta(t) - \eta(s) = \sum_{k \in \overline{E}} \theta_k d_k \tag{5.14}$$

and

$$\theta_k > 0 \quad \Longrightarrow \quad \varphi_k(u) = 0 \quad \text{for some } u \in (s, t]. \tag{5.15}$$

Substituting $\eta(t) - \eta(s)$ for $w$ in (3.2)–(3.3), and using the fact that $\mathcal{S} \subseteq \overline{E}$, the relation (5.14) implies that, for $i \in \mathcal{I}$,

$$\eta_i(t) - \eta_i(s) = \frac{\theta_i}{1 - \beta_i} - \beta_i \left( \sum_{j \in \overline{E}} \frac{\theta_j}{1 - \beta_j} \right)$$

$$= \frac{\theta_i}{1 - \beta_i} - \beta_i \left( \sum_{j \in \overline{E} \setminus \mathcal{S}} \frac{\theta_j}{1 - \beta_j} \right) - \beta_i \left( \sum_{j \in \mathcal{S}} \frac{\theta_j}{1 - \beta_j} \right). \tag{5.16}$$

Summing over $i \in \mathcal{S}$, this yields

$$\sum_{i \in \mathcal{S}} (\eta_i(t) - \eta_i(s)) = \left( \sum_{j \in \mathcal{S}} \frac{\theta_j}{1 - \beta_j} \right) \left( 1 - \sum_{i \in \mathcal{S}} \beta_i \right) - \left( \sum_{i \in \mathcal{S}} \beta_i \right) \left( \sum_{j \in \overline{E} \setminus \mathcal{S}} \frac{\theta_j}{1 - \beta_j} \right).$$

The last two relations can be combined with (5.12) to show that, for $i \in \mathcal{K}$,

$$\eta_i^{\mathcal{K}}(t) - \eta_i^{\mathcal{K}}(s) = \frac{\theta_i}{1 - \beta_i} - \beta_i \left( 1 + \frac{\sum_{k \in \mathcal{S}} \beta_k}{\sum_{k \in \mathcal{K}} \beta_k} \right) \left( \sum_{j \in \overline{E} \setminus \mathcal{S}} \frac{\theta_j}{1 - \beta_j} \right)$$

$$= \frac{\theta_i}{1 - \beta_i} - \beta_i^{\mathcal{K}} \left( \sum_{j \in \overline{E} \setminus \mathcal{S}} \frac{\theta_j}{1 - \beta_j} \right).$$

Using the definition

$$\theta_i^{\mathcal{K}} \doteq \frac{1 - \beta_i^{\mathcal{K}}}{1 - \beta_i} \theta_i \quad \text{for } i \in \mathcal{K} \tag{5.17}$$

and the fact that $\overline{E} \setminus \mathcal{S} = \overline{E} \cap \mathcal{K}$, the last equation can be rewritten as

$$\eta_i^{\mathcal{K}}(t) - \eta_i^{\mathcal{K}}(s) = \frac{\theta_i^{\mathcal{K}}}{1 - \beta_i^{\mathcal{K}}} - \beta_i^{\mathcal{K}} \left( \sum_{j \in \overline{E} \cap \mathcal{K}} \frac{\theta_j^{\mathcal{K}}}{1 - \beta_j^{\mathcal{K}}} \right) \quad \text{for } i \in \mathcal{K}.$$

From the definition of the $\mathcal{K}$-reduced directions of constraint given in (5.9), arguments analogous to those used to obtain (5.14) and the first line of (5.16) can be used to show that

$$\eta^{\mathcal{K}}(t) - \eta^{\mathcal{K}}(s) = \sum_{k \in \overline{E} \cap \mathcal{K}} \theta_k^{\mathcal{K}} d_k^{\mathcal{K}}.$$

That property 3 of the ESP is satisfied can now be inferred from the fact that, for $k \in \mathcal{K}$, (5.17) shows that $\theta_k^{\mathcal{K}} > 0$ if and only if $\theta_k > 0$, and (5.15) ensures that $\theta_k > 0$ implies $\varphi_k(u) = 0$ [equivalently, $\varphi_k^{\mathcal{K}}(u) = 0$] for some $u \in (s, t]$.

The fourth property of the ESP can be proved in an analogous fashion, splitting into the cases when $\varphi(t) = 0$ and $\varphi(t) \neq 0$ [in which case $H \neq d(\varphi(t))$ and the set $\{d_i : \varphi_i(t) = 0\}$ is linearly independent]. The argument is exactly analogous to the proof of property 3 and is thus omitted. Thus, we have shown that $\varphi^{\mathcal{K}} = \overline{\Gamma}^{\mathcal{K}}(\psi^{\mathcal{K}})$.

If $\psi \in \text{dom}(\Gamma)$, then, by Theorem 1.3(2) of [16], $\eta$ is of bounded variation on every compact time interval, which in turn implies that $\eta^{\mathcal{K}}$ is of bounded variation on every compact time interval. Since $(\varphi^{\mathcal{K}}, \eta^{\mathcal{K}})$ satisfy the $\mathcal{K}$-reduced ESP for $\psi^{\mathcal{K}}$, once again invoking Theorem 1.3(2) of [16], this shows that $(\varphi^{\mathcal{K}}, \eta^{\mathcal{K}})$ in fact solve the $\mathcal{K}$-reduced SP for $\psi^{\mathcal{K}}$, thus establishing the last statement of the theorem. $\square$

5.4. *Heavy traffic approximation for the unfinished work.* Fix $\nu \in \mathbb{R}^J$ with $\sum_{j \in \mathcal{I}} \nu_i = 0$ and denote the set of strictly subcritical classes $\mathcal{S}(\nu)$ simply by $\mathcal{S}$. Recall that, due to the ordering (4.16), $\mathcal{S}$ has the form $\{j_*, \ldots, J\}$ for some $j_* > 1$ (with $j_*$ set to $J+1$ if $\mathcal{S}$ is empty). Let $K = j_* - 1$, so that $K$ is the cardinality of the set $\mathcal{K} = \mathcal{I} \setminus \mathcal{S}$ of critical classes. It will be convenient to introduce the linear "projection" operators $\mathcal{L} : \mathcal{D}([0, \infty) : \mathbb{R}^J) \to \mathcal{D}([0, \infty) : \mathbb{R}^K)$ and $L : \mathbb{R}^J \to \mathbb{R}^K$ defined by

$$[\mathcal{L}f]_i = f_i + \beta_i^{\mathcal{K}}\left(\sum_{j \in \mathcal{S}} f_j\right) \qquad \text{if } i \in \mathcal{K},$$

for $f \in \mathcal{D}([0, \infty) : \mathbb{R}^J)$ and, analogously,

$$[Lv]_i = v_i + \beta_i^{\mathcal{K}}\left(\sum_{j \in \mathcal{S}} v_j\right) \qquad \text{if } i \in \mathcal{K}$$

for $v \in \mathbb{R}^J$. (In all cases, the sum over an empty set is taken to be equal to zero.) Define $\widehat{X}^{\mathcal{K}} \doteq \mathcal{L}(\hat{u} + B + \hat{c}\iota)$, and note that (since $\hat{u} \in \mathcal{M}$)

$$\widehat{X}_i^{\mathcal{K}} = \hat{u}_i + B_i + \hat{c}_i\iota + \beta_i^{\mathcal{K}} \sum_{j \in \mathcal{S}}(B_j + \hat{c}_j\iota).$$

Moreover, recalling the definition of the $\mathcal{K}$-reduced GPS ESM $\overline{\Gamma}^{\mathcal{K}}$ associated with the weight vector $\beta^{\mathcal{K}}$ defined in (5.8), let

$$\widehat{U}^{\mathcal{K}} \doteq \overline{\Gamma}^{\mathcal{K}}(\widehat{X}^{\mathcal{K}}),$$

and define the $J$-dimensional process $\widehat{U} \doteq (\widehat{U}^{\mathcal{K}}, 0^{\mathcal{S}})$, where here $0^{\mathcal{S}}$ is the identically zero process in the space $\Lambda_+^{\mathcal{S}} = \{(x_{j_*}, \ldots, x_J) : x \in \mathbb{R}^J\}$. In other words, let

(5.18) $$\widehat{U}_i \doteq \begin{cases} \widehat{U}_i^{\mathcal{K}}, & \text{for } i \in \mathcal{K}, \\ 0, & \text{for } i \in \mathcal{S}. \end{cases}$$

Also, define

$$\widehat{T} \doteq \hat{u} + B - \widehat{U} + \hat{c}\iota.$$

We now state and prove the main result of the paper.

THEOREM 5.7. *Suppose Assumptions* 2.1, 4.1, 5.1 *and the overall heavy traffic condition* (5.2) *are satisfied. Then* $\widehat{U}^n \Rightarrow \widehat{U}$ *and* $\widehat{T}^n \Rightarrow \widehat{T}$.

PROOF. For $n \in \mathbb{N}$, define $\mathcal{H}^{\mathcal{K},n} \doteq \mathcal{L}\mathcal{H}^n$, $\mathcal{X}^{\mathcal{K}n} \doteq \mathcal{L}\mathcal{X}^n$, $\mathcal{U}^{\mathcal{K},n} \doteq \Gamma^{\mathcal{K}}(\mathcal{X}^{\mathcal{K}n})$, $B^{\mathcal{K}} \doteq \mathcal{L}B$ and, likewise, $\alpha^{\mathcal{K}} \doteq L\alpha$, $\nu^{\mathcal{K}} \doteq L\nu$, $\hat{c}^{\mathcal{K}} \doteq L\hat{c}$, $\bar{u}^{\mathcal{K}} \doteq L\bar{u}$ and $\hat{u}^{\mathcal{K}} \doteq L\hat{u}$. Note that, since Assumption 4.1(1) and Assumption 5.1(1) together imply that $\bar{u} = 0$ and $\hat{u} \in \mathcal{M}$, it follows that $\bar{u}^{\mathcal{K}} = 0$ and $\hat{u}_i^{\mathcal{K}} = \hat{u}_i$ for $i \in \mathcal{K}$. By (4.20) and (4.21), we know that, for $i \in \mathcal{K}$, $\gamma_i - \alpha_i = \nu_i = -\beta_i \sum_{j \in \mathcal{S}} \nu_j / \sum_{k \in \mathcal{K}} \beta_k$. By (5.2), $\sum_{j=1}^J \nu_j = 0$. When rearranged, this shows that $\nu^{\mathcal{K}} = 0$ or, equivalently, that $\gamma^{\mathcal{K}} = \alpha^{\mathcal{K}}$. Also note that $\sum_{i=1}^K \alpha_i^{\mathcal{K}} = 1$. Now let $\overline{\mathcal{U}}^{\mathcal{K},n}$ and $\overline{\mathcal{X}}^{\mathcal{K},n}$ be the fluid scaled versions of $\mathcal{U}^{\mathcal{K},n}$ and $\mathcal{X}^{\mathcal{K},n}$, as defined in (4.1), and also define $\overline{\mathcal{X}}^{\mathcal{K}} \doteq 0$, $\hat{\mathcal{U}}^{\mathcal{K}} \doteq 0$. Then, by Assumption 4.1(1) and the definition of $\mathcal{U}_j^n(0)$ given after the proof of Theorem 5.4, we have

$$\lim_{n \to \infty} \overline{\mathcal{U}}^{\mathcal{K},n}(0) = L\left(\lim_{n \to \infty} \overline{\mathcal{U}}^n(0)\right) \le L\left(\lim_{n \to \infty} \frac{U^n(0)}{n}\right) = L\bar{u} = \bar{u}^{\mathcal{K}} = 0.$$

When combined with Theorem 5.4(2), Assumption 4.1(3), the linearity of the operators $\mathcal{L}$ and $L$ and the fact that $\nu^{\mathcal{K}} = 0$, this shows that $\overline{\mathcal{X}}^{\mathcal{K},n} \to \overline{\mathcal{X}}^{\mathcal{K}} \equiv 0$ u.o.c. Since $\Gamma^{\mathcal{K}}$ is the SM associated with the GPS model with basic and redistribution vectors $\alpha^{\mathcal{K}}$ and $\beta^{\mathcal{K}}$, Remark 2.4 ensures that it is Lipschitz continuous. Since $\overline{\mathcal{X}}^{\mathcal{K},n}$ and $\overline{\mathcal{X}}^{\mathcal{K}}$ are of bounded variation and thus lie in the domain of the SM, this allows us to conclude that $\overline{\mathcal{U}}^{\mathcal{K},n} \to \Gamma^{\mathcal{K}}(\overline{\mathcal{X}}^{\mathcal{K}}) = 0 \doteq \overline{\mathcal{U}}^{\mathcal{K}}$ u.o.c.

Now define $\widehat{\mathcal{X}}^{\mathcal{K},n}$, $\widehat{\mathcal{H}}^{\mathcal{K},n}$ and $\widehat{\mathcal{U}}^{\mathcal{K},n}$ as in (5.1), with $U, X$ and $H$ replaced by $\mathcal{U}^{\mathcal{K}}$, $\mathcal{X}^{\mathcal{K}}$ and $\mathcal{H}^{\mathcal{K}}$, respectively. The fact that $\overline{\mathcal{X}}^{\mathcal{K}} = \overline{\mathcal{U}}^{\mathcal{K}} = 0$, the scaling properties of the SM and the fact that the ESM coincides with the SM on the domain of the SP (see Remark 2.4) show that $\widehat{\mathcal{U}}^{\mathcal{K},n} = \overline{\Gamma}^{\mathcal{K}}(\widehat{\mathcal{X}}^{\mathcal{K},n})$ and

$$\widehat{\mathcal{X}}^{\mathcal{K},n} = \frac{L\mathcal{U}^n(0)}{\sqrt{n}} + \mathcal{L}(\widehat{\mathcal{H}}^n + \sqrt{n}[\nu + (\gamma^n - \gamma)]\iota).$$

Along with Theorem 5.4(3), Assumption 5.1(1) and Assumption 5.1(3), this shows that $\widehat{\mathcal{X}}^{\mathcal{K},n} \to \widehat{X}^{\mathcal{K}}$ u.o.c. The fact that $\widehat{\mathcal{U}}^{\mathcal{K},n} = \overline{\Gamma}^{\mathcal{K}}(\widehat{\mathcal{X}}^{\mathcal{K},n})$ and the Lipschitz continuity of the GPS ESM $\overline{\Gamma}^{\mathcal{K}}$ (which follows from Remark 2.4) then shows that $\widehat{\mathcal{U}}^{\mathcal{K},n} \to \widehat{U}^{\mathcal{K}}$ u.o.c. Now fix $\omega \in \Omega$. Since Lemma 5.5, along with Theorem 5.6, shows that there exists $N = N(\omega) < \infty$ such that, for all $n \geq N(\omega)$,

$$\widehat{\mathcal{U}}_i^n(\omega) = \begin{cases} 0, & \text{if } i \in \mathcal{S}, \\ [\Gamma^{\mathcal{K}}(\widehat{\mathcal{X}}^{\mathcal{K},n}(\omega))]_i = \widehat{\mathcal{U}}_i^{\mathcal{K},n}(\omega), & \text{if } i \in \mathcal{K}, \end{cases}$$

we conclude that $\widehat{\mathcal{U}}^n(\omega) \to \widehat{U}(\omega)$ for every $\omega \in \Omega$.

Furthermore, since $\overline{\mathcal{U}} = \widehat{U} \equiv 0$ and $\Gamma$ is Lipschitz continuous, there exists a constant $C < \infty$ such that, for any $T < \infty$,

$$\sup_{t \in [0,T]} |\widehat{U}^n(t) - \widehat{\mathcal{U}}^n(t)|$$

$$= \sqrt{n} \sup_{t \in [0,T]} |\overline{U}^n(t) - \overline{\mathcal{U}}^n(t)|$$

$$= \sqrt{n} \sup_{t \in [0,T]} |\Gamma(\overline{X}^n)(t) - \Gamma(\overline{\mathcal{X}}^n)(t)|$$

$$\leq C\sqrt{n} \sup_{t \in [0,T]} |\overline{X}^n(t) - \overline{\mathcal{X}}^n(t)|$$

$$\leq C \left| \frac{U^n(0)}{\sqrt{n}} - \frac{\mathcal{U}^n(0)}{\sqrt{n}} \right| + C\sqrt{n} \sup_{t \in [0,T]} |\overline{H}^n(t) - \overline{\mathcal{H}}^n(t)|$$

$$\leq JC \max_{j \in \mathcal{S}} \frac{U_j^n(0)}{\sqrt{n}} + C \sup_{t \in [0,T]} |\widehat{H}^n(t) - \widehat{\mathcal{H}}^n(t)|.$$

Taking limits as $n \to \infty$, Assumption 5.1(1) and the fact that $\hat{u} \in \mathcal{M}$ show that the first term in the last line of the display tends to zero, while Assumption 5.1(2) and Theorem 5.4(3) show that the second term tends to zero. When combined with the fact that $\widehat{\mathcal{U}}^n \to \widehat{U}$ u.o.c., this shows that $\widehat{U}^n \to \widehat{U}$ u.o.c., as desired.

Finally, using equations (2.1) and (2.2), for $n \in \mathbb{N}$ and $t \in [0, \infty)$, we can write

$$T_i^n(t) = U_i^n(0) + H_i^n(t) - U_i^n(t).$$

Using the relation $\overline{T} = \gamma \iota$, which follows from Theorem 3.2, this implies that

$$\overline{T}_i^n(t) = \frac{U_i^n(0)}{\sqrt{n}} + \overline{H}_i^n(t) - \overline{U}_i^n(t) + \sqrt{n}(\gamma_i^n - \gamma_i)t.$$

By Assumption 5.1 and the diffusion limit for the unfinished work just proved above, we then have

(5.19) $$\widehat{T}^n \to \widehat{T} \doteq \hat{u} + B - \widehat{U} + \hat{c}\iota,$$

u.o.c. By Remark 5.2, the pathwise u.o.c. limits must be replaced by weak convergence, and the proof of the theorem is complete. □

REMARK 5.8. It is also possible to introduce additional primitives to describe the queue length $Q^n$, sojourn time $V^n$ and waiting time $W^n$ processes associated with this model. Indeed, this was done in [17] where, under general assumptions on the first-order properties of these primitives (see Assumption 4.5 and Condition 4.8 of [17]), fluid limits for these processes were established (see Theorems 4.7 and 4.11 of [17]). In particular, it was also shown in Lemma 4.13 of [17] that when the heavy traffic condition $\sum_{i=1}^{J} \gamma_i = 1$ holds (and the initial conditions converge in a suitable manner to 0), the fluid limits of all these processes are identically zero. In addition, under general functional central limit type assumptions (see Assumptions 4.12 and 4.16 of [17]), diffusion limits for these processes were deduced from the corresponding diffusion limits for the unfinished work and busy time processes. Indeed, it was shown in Theorems 4.18 and 4.19 of [17] that the diffusion limit $\widehat{Q}$ for the queue length process satisfies $\widehat{Q}_i = \mu_i \widehat{U}_i$ for $i \in \mathcal{I}$, where $1/\mu_i$ corresponds to the long-run average service requirement of a class $i$ customer, while the diffusion limits $\widehat{V}$ and $\widehat{W}$ for the sojourn and waiting time processes, respectively, satisfy $\widehat{V}_i = \widehat{W}_i = \widehat{U}_i/\gamma_i$ for $i \in \mathcal{I}$. These diffusion results were obtained in [17] under the balanced heavy traffic condition $\gamma = \alpha$. However, it can be shown that these results continue to hold even in the unbalanced case considered in this paper. Indeed, they can be deduced from the results of Theorem 5.7 by following almost verbatim the proofs of Theorems 4.18 and 4.19 in [17], with the only modification that $\alpha$ be replaced by $\gamma$ in those proofs. As a consequence, we omit a detailed exposition of these results.

**Acknowledgments.** The authors are grateful to comments by the referees that led to an improvement in the exposition of the paper.

Department of Mathematical Sciences  
Carnegie Mellon University  
Pittsburgh, Pennsylvania 15213  
USA  
E-mail: kramanan@math.cmu.edu

Alcatel–Lucent Bell Labs  
600 Mountain Avenue  
Murray Hill, New Jersey 07974  
USA  
E-mail: marty@alcatel-lucent.com